\documentclass[11pt]{article}
\usepackage[tbtags]{amsmath}
\usepackage{amssymb}
\usepackage{amsthm}
\usepackage[misc]{ifsym}
\usepackage{cases}
\usepackage{mathrsfs}
\numberwithin{equation}{section}
\normalsize
\title{\bf Linear Quadratic Stackelberg Stochastic Differential Games: Closed-Loop Solvability
\thanks{This work is financially supported by the National Key R\&D Program of China (2018YFB1305400), and the National Natural Science Foundations of China (11971266, 11571205, 11831010).}}
\author{\normalsize Zixuan Li\thanks{\it School of Mathematics, Shandong University, Jinan 250100, P.R. China, E-mail: 201812064@mail.sdu.edu.cn},\quad
Jingtao Shi\thanks{Corresponding author,\ \it School of Mathematics, Shandong University, Jinan 250100, P.R. China, E-mail: shijingtao@sdu.edu.cn}}
\newtheorem{mypro}{Proposition}[section]
\newtheorem{mythm}{Theorem}[section]
\newtheorem{mydef}{Definition}[section]

\newtheorem{Remark}{Remark}[section]
\begin{document}
\maketitle

\noindent{\bf Abstract:}\quad This paper is concerned with the closed-loop solvability of one kind of linear-quadratic Stackelberg stochastic differential game, where the coefficients are deterministic. The notion of the closed-loop solvability is introduced, which require to be independent of the initial state. The follower's problem is solved first, and the closed-loop optimal strategy is characterized by a Riccati equation, together with an adapted solution to a linear backward stochastic differential equation. Then the necessary conditions of the existence of the leader's nonanticipating closed-loop optimal strategy is obtained via a system of cross-coupled Riccati equations. The sufficiency is open since the completion-of-square method is invalid.

\vspace{2mm}


\vspace{2mm}


\section{Introduction}

Let us first introduce some notations which will be used throughout the paper.

Let $T>0$ be a finite time duration. Let $\mathbb{R}^{n \times m}$ be the set of all $(n \times m)$ matrices, $\mathbb{S}^n$ be the set of all $(n \times n)$ symmetric matrices. For a Banach space $H$ (for example, $H=\mathbb{R}^n,\mathbb{R}^{n \times m},\mathbb{S}^n$), let $L^p(0,T;H)\,(1 \leqslant p \leqslant \infty)$ be the space of all $H$-valued functions which are $L^p$-integrable on $[0,T]$, and $C([0,T];H)$ be the space of all $H$-valued continuous functions on $[0,T]$.

Let $(\Omega,\mathcal{F},\mathbb{F},\mathbb{P})$ be a completed filtered probability space on which a standard one-dimensional Brownian motion $W=\{W(t);0 \leqslant t < \infty \}$ is defined, where $\mathbb{F}=\{\mathcal{F}_t\}_{t\geqslant 0}$ is natural filtration of $W$ augmented by all the $\mathbb{P}$-null sets in $\mathcal{F}$. We denote
\begin{equation*}
\begin{aligned}
&L^2_{\mathcal{F}_t}(\Omega;H)=\Big\{\xi:\Omega \to H\,|\,\xi\,\,\mbox{is}\,\,\mathcal{F}_t\mbox{-measurable},\,\,\,\mathbb{E}|\xi|^2<\infty \Big\},\,\, t\in(0,T],\\
&L^2_{\mathbb{F}}(0,T;H)=\bigg\{X(\cdot):[0,T]\times\Omega\to H\,\big|\,X(\cdot)\,\,\mbox{is}\,\,\mathbb{F}\mbox{-progressively measurable},\\
&\hspace{3.5cm}\,\,\,\mathbb{E}\int_0^T|X(s)|^2ds<\infty \bigg\}.
\end{aligned}
\end{equation*}

We consider the following controlled linear {\it stochastic differential equation} (SDE for short):
\begin{equation}\label{state}\left\{
\begin{aligned}
dx^{u_1,u_2}(s)&=\big[A(s)x^{u_1,u_2}(s)+B_1(s)u_1(s)+B_2(s)u_2(s)\big]ds\\
               &\quad+C(s)x^{u_1,u_2}(s)dW(s),\ s \in [0,T],\\
 x^{u_1,u_2}(0)&=x,
\end{aligned}
\right.\end{equation}
where $x\in\mathbb{R}^n$ is the given initial state, $A(\cdot),B_1(\cdot),B_2(\cdot),C(\cdot)$ are given deterministic matrix-valued functions of proper dimensions. In the above, $x^{u_1,u_2}(\cdot)$ is the \textit{state process} with values in $\mathbb{R}^n$, and $u_1(\cdot),u_2(\cdot)$ are \textit{control processes} with values in $\mathbb{R}^{m_1}$ and $\mathbb{R}^{m_2}$, taken by the two players in the games, labeled 1 and 2, respectively. We introcuce the following Hilbert space:
\begin{equation}\label{ac}
\begin{aligned}
\mathcal{U}_i[0,T]=\bigg\{u_i:\,&[0,T] \times \Omega \to \mathbb{R}^{m_i} \,\,\big|\,\,u_i(\cdot)\,\, \mbox{is}\, \,\mathbb{F}\mbox{-progressively measurable}, \\
                                &\mathbb{E}\int_{0}^{T} |u_i(s)|^2 ds < \infty \bigg\},\,\quad \,i=1,2.
\end{aligned}
\end{equation}
The control processes $u_1(\cdot)\in\mathcal{U}_1[0,T]$ and $u_2(\cdot)\in\mathcal{U}_2[0,T]$ are called \textit{adimissible controls}.
\par Under some mild conditions on the coefficients, for any $(x,u_1(\cdot),u_2(\cdot)) \in \mathbb{R}^n \times \mathcal{U}_1[0,T] \times \mathcal{U}_2[0,T]$, there exists a unique (strong) solution $x^{u_1,u_2}(\cdot)\in L^2_{\mathcal{F}}(0,T;\mathbb{R}^n)$ to (\ref{state}). Thus, we can define the \textit{cost functionals} for the players as follows: For $i=1,2$,
\begin{equation}\label{cfi}
\begin{aligned}
J_i(x;u_1(\cdot),u_2(\cdot))&=\mathbb{E}\bigg\{ \int_0^T \Big[\big\langle Q_i(s)x^{u_1,u_2}(s),x^{u_1,u_2}(s)\big\rangle+\big\langle R_i(s)u_i(s),u_i(s)\big\rangle\Big] dt\\
                            &\qquad\qquad\quad +\big\langle G_ix^{u_1,u_2}(T),x^{u_1,u_2}(T)\big\rangle \bigg\},
\end{aligned}
\end{equation}
where $Q_i(\cdot),R_i(\cdot)$ are deterministic matrix-valued functions of proper dimensions with $Q_i(\cdot)^\top =Q_i(\cdot),\,\,R_i(\cdot)^\top =R_i(\cdot)$, and $G_i$ is a symmetric matrix.

\par In our Stackelberg game framework, Player 1 is the follower and Player 2 is the leader. For any choice $u_2(\cdot) \in \mathcal{U}_2[0,T]$ of Player 2 and a fixed initial state $x \in \mathbb{R}^n$, Player 1 would like to choose a $\bar{u}_1(\cdot) \in \mathcal{U}_1[0,T]$ so that $J_1(x;\bar{u}_1(\cdot),u_2(\cdot))$ is the minimum of $J_1(x;u_1(\cdot),u_2(\cdot))$ over $u_1(\cdot) \in \mathcal{U}_1[0,T]$. Knowing Player 1 would take such an optimal control $\bar{u}_1(\cdot)$, Player 2 would like to choose a $\bar{u}_2(\cdot) \in \mathcal{U}_2[0,T]$ to minimize $J_2(x;\bar{u}_1(\cdot),u_2(\cdot))$ over $u_2(\cdot) \in \mathcal{U}_2[0,T]$. We refer to such a problem as a {\it linear quadratic (LQ for short) Stackelberg stochastic differential game}.

In a more rigorous way, Player 1 wants to find a map $\bar{u}_1:\mathcal{U}_2[0,T] \times \mathbb{R}^n \to \mathcal{U}_1[0,T] $ and Player 2 want to find a $\bar{u}_2(\cdot) \in \mathcal{U}_2[0,T]$ such that
\begin{equation}\label{problem}
\left\{
\begin{aligned}
            J_1(x;\bar{u}_1[u_2,x](\cdot),u_2(\cdot))&=\mathop{\min}\limits_{u_1(\cdot) \in \mathcal{U}_1[0,T]}J_1(x;u_1(\cdot),u_2(\cdot)),\quad \,\,\forall u_2(\cdot) \in \mathcal{U}_2[0,T],\\
J_2(x;\bar{u}_1[\bar{u}_2,x](\cdot),\bar{u}_2(\cdot))&=\mathop{\min}\limits_{u_2(\cdot) \in \mathcal{U}_2[0,T]}J_2(x;\bar{u}_1[u_2,x](\cdot),u_2(\cdot)).
\end{aligned}
\right.
\end{equation}
If the above pair $(\bar{u}_1[\cdot,x](\cdot),\bar{u}_2(\cdot))$ exists, we refer to it as an \textit{open-loop solution} to the above LQ Stackelberg stochastic differential game, for $x\in\mathbb{R}^n$.

\par The theory of Stackelberg game can be traced back to Stackelberg \cite{Stackelberg1934}, who put forward the Stackelberg game and the concept of Stackelberg solution in static competitive economics with a hierarchical structure. Simann and Cruz \cite{SC1973} studied the multi-stages and dynamic LQ Stackelberg differential games, where feedback Stackelberg solutions are introduced. Castanon and Athans \cite{CA1976} considered feedback Stackelberg strategies for two person linear multistage games with quadratic performance criteria and noisy measurements and gave a explicit solutions when the information sets are nested in a stochastic case. Ba\c{s}ar and Selbuz \cite{BS1979} considered the closed-loop Stackelberg solution to a class of LQ two-person nonzero sum differential games. Bagchi and Ba\c{s}ar \cite{BB1981} investigated the LQ stochastic Stackelberg differential game, where the diffusion term of the state equation does not contain the state and control variables. Yong \cite{Yong2002} extended the LQ stochastic Stackelberg differential game to random and state-control dependent coefficients, and obtained the feedback representation of the open-loop solution via some stochastic Riccati equations. In the past decades, there have been a great deal of works on this issue, for jump diffusions see \O ksendal et al. \cite{OSU2013}, Moon \cite{Moon2020}, for different information structures see Ba\c{s}ar and Olsder \cite{BO1998}, Bensoussan et al. \cite{BCS2015}, for time-delayed systems see Xu and Zhang \cite{XZ2016}, Xu et al. \cite{XSZ2018}, for mean field's type models related with multiple followers and large populations see Mukaidani and Xu \cite{MX2015}, Moon and Ba\c{s}ar \cite{MB2018}, Li and Yu \cite{LiYu2018}, Lin et al. \cite{LJZ2019}, Wang and Zhang \cite{WZ2020}, for partial/asymmetric/overlapping information see Shi et al. \cite{SWX2016,SWX2017,SWX2020}, for backward stochastic systems see Du and Wu \cite{DW2019}, Zheng and Shi \cite{ZS2020}, for time-inconsistent case see Moon and Yang \cite{MY2021}.

Our interest in this paper lies in the {\it closed-loop solution} or the {\it closed-loop solvability} for the above LQ Stackelberg stochastic differential game. To our best knowledge, this topic has not been studied in the literature yet. However, the closed-loop solution for (LQ) Stackelberg stochastic differential game are mentioned but not addressed in \cite{Yong2002} and \cite{BCS2015}. In 2014, Sun and Yong \cite{SY2014} introduce the notions of open-loop and closed-loop solvabilities for an LQ stochastic optimal control problem, which is a special case when only one player/controller is considered for open-loop and closed-loop saddle points for an LQ two-person zero-sum stochastic differential game. Sun et al. \cite{SLY2016} further gives more detailed characterizations of the closed-loop solvability for the LQ stochastic optimal control problem. Sun and Yong \cite{SY2019} is devoted to the open-loop and closed-loop Nash equilibria for an LQ two-person nonzero-sum stochastic differential game. The existence of an optimal closed-loop strategy for an LQ mean-field optimal control problem is studied in Li et al. \cite{LSY2016}. Sun and Yong \cite{SY2018} obtained the equivalence of open-loop and closed-loop solvabilities for the LQ stochastic optimal control problem in an infinite horizon. Very recently, Li et al. \cite{LSY2020} extended the previous results to LQ mean-field two-person zero-sum and nonzero sum stochastic differential games in an infinite horizon.

In this paper, we first solve the follower's problem, and his closed-loop optimal strategy is characterized by a Riccati equation, together with the adapted solution to a linear {\it backward stochastic differential equation} (BSDE for short). Then we solve the leader's problem, whose state equation is a {\it forward-backward stochastic differential equation} (FBSDE for short). We introduce the definition of the closed-loop solvability of the leader's problem. Necessary conditions for the nonanticipating closed-loop optimal strategy of the leader are given, via a cross-coupled Riccati equation system. The sufficiency is open since the completion-of-square method is invalid.

The rest of this paper is organized as follows. Section 2 gives some preliminaries, to introduce the closed-loop solution to the LQ Stackelberg stochastic differential game. Section 3 is devoted to solve the problem of the follower. The sufficient and necessary conditions for the closed-loop solvability of the follower's problem are given. In Section 4, necessary conditions for the closed-loop solvability of the leader's problem is given. Finally, in Section 5 some concluding remarks are given.

%Problem(SLQ)
\section{Preliminaries}

First of all, we recall the open-loop and closed-loop solvabilities for the LQ stochastic optimal control problem (see \cite{SLY2016}, for example). Consider the linear state equation
\begin{equation}\label{SLQ}
\left\{\begin{aligned}
dX^u(s)&=\big[A(s)X^u(s)+B(s)u(s)+b(s)\big]ds\\
       &\quad +\big[C(s)X^u(s)+D(s)u(s)+\sigma(s)\big]dW(s),\ \ s \in [0,T],\\
 X^u(0)&=x,
\end{aligned}\right.
\end{equation}
and the quadratic cost functional:
\begin{equation}\label{cf}
\begin{split}
J(x;u(\cdot))&=\mathbb{E} \bigg\{\big\langle GX^u(T),X^u(T)\big\rangle\\
&\qquad+\int_0^T \bigg[\bigg\langle
\left( \begin{array}{cc} Q(s) & S(s)^\top \\ S(s) & R(s)\end{array} \right)
\left( \begin{array}{c} X^u(s) \\ u(s)\end{array} \right),
\left( \begin{array}{c} X^u(s) \\ u(s)\end{array}\right)\bigg\rangle \bigg]ds \bigg\}.
\end{split}
\end{equation}
We adopt the following assumptions.
\par \textbf{(S1)} The coefficients of the state equation (\ref{SLQ}) satisfy the following:
\begin{equation}\nonumber
\begin{cases}
A(\cdot) \in L^1(0,T;\mathbb{R}^{n \times n}),\,\,B(\cdot) \in L^2(0,T;\mathbb{R}^{n \times m}),\,\,b(\cdot) \in L^2_\mathbb{F}(\Omega;L^1(0,T;\mathbb{R}^n)),\\
C(\cdot) \in L^2(0,T;\mathbb{R}^{n \times n}),\,\,D(\cdot) \in L^\infty(0,T;\mathbb{R}^{n \times m}),\,\,\sigma(\cdot) \in L^2_\mathbb{F}(0,T;\mathbb{R}^n)\\
\end{cases}
\end{equation}
\par \textbf{(S2)} The weighting coefficients of the cost functional (\ref{cf}) satisfy the following:
\begin{equation}\nonumber
Q(\cdot) \in L^1(0,T;\mathbb{S}^n),\,\,S(\cdot) \in L^2(0,T;\mathbb{R}^{m \times n}),\,\,\,R(\cdot) \in L^{\infty}(0,T;\mathbb{S}^{m}),\,\,\,G \in \mathbb{S}^n.\\
\end{equation}

Under (S1) and (S2), for any $x \in \mathbb{R}^n$ and $u(\cdot) \in \mathcal{U}[0,T] \equiv L^2_{\mathbb{F}}(0,T;\mathbb{R}^m)$, the state equation (\ref{SLQ}) admits a unique strong solution $X^u(\cdot)\in L^2_\mathbb{F}(0,T;\mathbb{R}^n)$ and the cost functional (\ref{cf}) is well-defined. Therefore, the following problem is meaningful.

\textbf{Problem (SLQ)}. For any initial state $x \in \mathbb{R}^n$, find a $\bar{u}(\cdot) \in \mathcal{U}[0,T]$ such that
\begin{equation}\label{PLQ}
J(x;\bar{u}(\cdot))=\mathop{\inf}\limits_{u(\cdot) \in \mathcal{U}[0,T]}J(x;u(\cdot)) \equiv V(x).
\end{equation}
Any $\bar{u}(\cdot) \in \mathcal{U}[0,T]$ satisfying (\ref{PLQ}) is called \textit{an open-loop optimal control} of Problem (SLQ) for $x$, the corresponding $\bar{X}(\cdot)\equiv X^{\bar{u}}(\cdot)$ is called \textit{an open-loop optimal state process} and $(\bar{X}(\cdot),\bar{u}(\cdot))$ is called \textit{an open-loop optimal pair}. $V(\cdot)$ is called \textit{the value function} of Problem (SLQ).

%SLQ
\begin{mydef}\label{def2.1}
Let $x\in \mathbb{R}^n$. If there exists a (unique) $\bar{u}(\cdot) \in \mathcal{U}[0,T]$ such that (\ref{PLQ}) holds, then we say that Problem (SLQ) is (uniquely) open-loop solvable at $x$. If Problem (SLQ) is (uniquely) open-loop solvable for every $x \in \mathbb{R}^n$, then we say that Problem (SLQ) is (uniquely) open-loop solvable.
\end{mydef}

The following result is concerned with open-loop solvability of Problem (SLQ), whose proof can be found in \cite{SLY2016} (see also \cite{SY2014}).

% SLQ open-loop optimal solvability
\begin{mypro}\label{opsn}
Let (S1)-(S2) hold. For an initial state $x \in \mathbb{R}^n$, a state-control pair $(\bar{X}(\cdot),\bar{u}(\cdot))$ is an open-loop optimal pair of Problem (SLQ) if and only if the following hold:\\
(\romannumeral 1) The stationarity condition holds:
\begin{equation}\label{sc}
B(s)^\top\bar{Y}(s)+D(s)^\top\bar{Z}(s)+S(s)\bar{X}(s)+R(s)\bar{u}(s)=0,\quad a.e.\ s \in[0,T],\ \mathbb{P}\mbox{-}a.s.,
\end{equation}
where $(\bar{Y}(\cdot),\bar{Z}(\cdot)) \in L^2_{\mathbb{F}}(0,T;\mathbb{R}^n) \times L^2_{\mathbb{F}}(0,T;\mathbb{R}^n)$ is the solution to the following BSDE:
\begin{equation}\left\{
\begin{aligned}
d\bar{Y}(s)&=-\big[A(s)^\top \bar{Y}(s)+C(s)^\top\bar{Z}(s)+Q(s)\bar{X}(s)+S(s)^\top\bar{u}(s)\big]ds\\
           &\quad+\bar{Z}(s)dW(s),\,\,\,\, s \in [0,T],\\
 \bar{Y}(T)&=G\bar{X}(T).
\end{aligned}\right.
\end{equation}
(\romannumeral 2) The map $u(\cdot) \to J(0;u(\cdot))$ is convex.
\end{mypro}

Next, take $\Theta(\cdot) \in L^2(0,T;\mathbb{R}^{m \times n}) \equiv \mathcal{Q}[0,T]$ and $v(\cdot) \in \mathcal{U}[0,T]$. For any $x \in \mathbb{R}^n$, let us consider the following linear equation:
\begin{equation}\label{cl}
\left\{
\begin{aligned}
dX^{\Theta,v}(s)&=\Big\{\big[A(s)+B(s)\Theta(s)\big]X^{\Theta,v}(s)+B(s)v(s)+b(s)\Big\}ds\\
                &\quad +\Big\{\big[C(s)+D(s)\Theta(s)\big]X^{\Theta,v}(s)+D(s)v(s)+\sigma(s)\Big\}dW(s),\quad s \in [0,T],\\
 X^{\Theta,v}(0)&=x,
\end{aligned}\right.
\end{equation}
which admits a unique solution $X^{\Theta,v}(\cdot)\in L^2_\mathbb{F}(0,T;\mathbb{R}^n)$, depending on the $\Theta(\cdot)$ and $v(\cdot)$. The above equation (\ref{cl}) is called a {\it closed-loop system} of the original state equation (\ref{SLQ}) under a {\it closed-loop strategy} $(\Theta(\cdot),v(\cdot))$. We point out that $(\Theta(\cdot),v(\cdot))$ is independent of the initial state $x \in \mathbb{R}^n$. With the above solution $X^{\Theta,v}(\cdot)$, we define
\begin{equation}\label{ccf}
\begin{aligned}
&J(x;\Theta(\cdot)X^{\Theta,v}(\cdot)+v(\cdot))=\mathbb{E} \bigg\{ \big\langle GX^{\Theta,v}(T),X^{\Theta,v}(T)\big\rangle \\
&+\int_0^T \bigg[\bigg\langle
\left( \begin{array}{cc} Q(s) & S(s)^\top \\ S(s) & R(s)\end{array} \right)
\left( \begin{array}{c} X^{\Theta,v}(s) \\ \Theta(s)X^{\Theta,v}(s)+v(s)\end{array} \right),
\left( \begin{array}{c} X^{\Theta,v}(s) \\ \Theta(s)X^{\Theta,v}(s)+v(s)\end{array}\right)\bigg\rangle \bigg]ds \bigg\},
\end{aligned}
\end{equation}
and recall the following definition.
\begin{mydef}\label{def2.2}
A pair $(\bar{\Theta}(\cdot),\bar{v}(\cdot)) \in \mathcal{Q}[0,T] \times \mathcal{U}[0,T]$ is called a \textit{closed-loop optimal strategy} of Problem (SLQ) if
\begin{equation}\label{CLQ}
\begin{split}
J(x;\bar{\Theta}(\cdot)X^{\bar{\Theta},\bar{v}}(\cdot)+\bar{v}(\cdot)) \leqslant J(x;\Theta(\cdot)X^{\Theta,v}(\cdot)+v(\cdot)),\\
\forall x \in \mathbb{R}^n,\,\,\forall  (\Theta(\cdot),v(\cdot)) \in \mathcal{Q}[0,T] \times \mathcal{U}[0,T].
\end{split}
\end{equation}
If there exists a (unique) pair $(\bar{\Theta}(\cdot),\bar{v}(\cdot)) \in \mathcal{Q}[0,T] \times \mathcal{U}[0,T]$ such that (\ref{CLQ}) holds, we say that Problem (SLQ) is (uniquely) closed-loop solvable.
\end{mydef}

We emphasize that the pair $(\bar{\Theta}(\cdot),\bar{v}(\cdot))$ is required to be independent of the initial state $x \in \mathbb{R}^n$. The following results about some equivalent definitions is also from \cite{SLY2016}.

\begin{mypro}\label{relation}
Let (S1)-(S2) hold and $(\bar{\Theta}(\cdot),\bar{v}(\cdot)) \in \mathcal{Q}[0,T] \times \mathcal{U}[0,T]$. Then the following statements are equivalent:\\
(\romannumeral 1) $(\bar{\Theta}(\cdot),\bar{v}(\cdot))$ is a closed-loop optimal strategy of Problem (SLQ).\\
(\romannumeral 2) For any $x \in \mathbb{R}^n$ and $v(\cdot) \in \mathcal{U}[0,T]$,
\begin{equation}\label{equlvalent}
J(x;\bar{\Theta}(\cdot)X^{\bar{\Theta},\bar{v}}(\cdot)+\bar{v}(\cdot)) \leqslant J(x;\bar{\Theta}(\cdot)X^{\bar{\Theta},v}(\cdot)+v(\cdot)).\\
\end{equation}
(\romannumeral 3) For any $x \in \mathbb{R}^n$ and $u(\cdot) \in \mathcal{U}[0,T]$,
\begin{equation}\label{open and close relation}
J(x;\bar{\Theta}(\cdot)X^{\bar{\Theta},\bar{v}}(\cdot)+\bar{v}(\cdot)) \leqslant J(x;u(\cdot)).
\end{equation}	
\end{mypro}

From the above result, we see that if $(\bar{\Theta}(\cdot),\bar{v}(\cdot))$ is a closed-loop optimal strategy of Problem (SLQ), then for any fixed initial state $x \in \mathbb{R}^n$, (\ref{open and close relation}) implies that the {\it outcome}
\begin{equation}\nonumber
\bar{u}(\cdot)=\bar{\Theta}(\cdot)X^{\bar{\Theta},\bar{v}}(\cdot)+\bar{v}(\cdot) \in \mathcal{U}[0,T]
\end{equation}
is an open-loop optimal control of Problem (SLQ) for $x$. Therefore, for Problem (SLQ), the closed-loop solvability implies the open-loop solvability for any $x \in \mathbb{R}^n$.

\vspace{1mm}

We now return to our LQ Stackelberg stochastic differential game (\ref{state})--(\ref{problem}). We denote $L^2(0,T;\mathbb{R}^{m_i \times n}) \equiv \mathcal{Q}_i[0,T]$ for $i=1,2$.

First, for any $u_2(\cdot)\in \mathcal{U}_2[0,T]$, take $\Theta_1(\cdot) \in \mathcal{Q}_1[0,T]$ and $v_1(\cdot)\in \mathcal{U}_1[0,T]$. For any $x \in \mathbb{R}^n$, let us consider the following linear equation:
\begin{equation}\label{cse}\left\{
\begin{aligned}
dx^{\Theta_1,v_1,u_2}(s)&=\Big\{\big[A(s)+B_1(s)\Theta_1(s)\big]x^{\Theta_1,v_1,u_2}(s)+B_1(s)v_1(s)+B_2(s)u_2(s)\Big\}ds\\
                        &\quad+C(s)x^{\Theta_1,v_1,u_2}(s)dW(s),\quad s \in [0,T],\\
 x^{\Theta_1,v_1,u_2}(0)&=x,
\end{aligned}\right.
\end{equation}
which admits a unique solution $x^{\Theta_1,v_1,u_2}(\cdot)\in L^2_\mathbb{F}(0,T;\mathbb{R}^n)$, depending on the $\Theta_1(\cdot)$ and $v_1(\cdot)$. The above is called a {\it closed-loop system} of the original state equation (\ref{state}) under the {\it closed-loop strategy} $(\Theta_1(\cdot),v_1(\cdot))$ of the follower. We point out that $(\Theta_1(\cdot),v_1(\cdot))$ is independent of the initial state $x$. With the above solution $x^{\Theta_1,v_1,u_2}(\cdot)$, we define
\begin{equation}\label{fccf}
\begin{aligned}
&J_1(x;\Theta_1(\cdot)x^{\Theta_1,v_1,u_2}(\cdot)+v_1(\cdot),u_2(\cdot))=\mathbb{E} \bigg\{ \big\langle G_1x^{\Theta_1,v_1,u_2}(T),x^{\Theta_1,v_1,u_2}(T)\big\rangle\\
&\qquad+\int_0^T\Big[ \big\langle \big[Q_1(s)+\Theta_1^\top(s) R_1(s)\Theta_1(s)\big]x^{\Theta_1,v_1,u_2}(s),x^{\Theta_1,v_1,u_2}(s) \big\rangle\\
&\qquad+2\big\langle R_1(s)\Theta_1(s)x^{\Theta_1,v_1,u_2}(s),v_1(s)\big\rangle+\big\langle R_1(s)v_1(s),v_1(s)\big\rangle ds\bigg\}
\end{aligned}
\end{equation}and introduce the following notion.
\begin{mydef}\label{def2.3}
A quadruple $(\bar{\Theta}_1(\cdot),\bar{v}_1(\cdot),\bar{\Theta}_2(\cdot),\bar{v}_2(\cdot)) \in \mathcal{Q}_1[0,T] \times \mathcal{U}_1[0,T] \times \mathcal{Q}_2[0,T] \times \mathcal{U}_2[0,T] $ is called a (unique) closed-loop solution to our LQ Stackelberg stochastic differential game, if\\
(\romannumeral 1) For any $x \in \mathbb{R}^n$ and given $u_2(\cdot) \in \mathcal{U}_2[0,T]$, Player 1 could find two maps: $\bar{\Theta}_1:\mathcal{U}_2[0,T]\rightarrow\mathcal{Q}_1[0,T]$ and $\bar{v}_1:\mathcal{U}_2[0,T]\rightarrow\mathcal{U}_1[0,T]$ such that
\begin{equation}
\begin{aligned}
&J_1(x;\bar{\Theta}_1[u_2](\cdot)\bar{x}^{u_2}(\cdot)+\bar{v}_1[u_2](\cdot),u_2(\cdot))\leqslant J_1(x;\Theta_1[u_2](\cdot)x^{\Theta_1[u_2],v_1[u_2],u_2}(\cdot)+v_1[u_2](\cdot),u_2(\cdot)),\\
&\hspace{6cm} \forall\, \Theta_1:\mathcal{U}_2[0,T]\rightarrow\mathcal{Q}_1[0,T],\,\,v_1:\mathcal{U}_2[0,T]\rightarrow\mathcal{U}_1[0,T],\\
\end{aligned}
\end{equation}
where $\bar{x}^{u_2}(\cdot)\equiv x^{\bar{\Theta}_1[u_2],\bar{v}_1[u_2],u_2}(\cdot)$.\\
(\romannumeral 2) There exist a (unique) pair $(\bar{\Theta}_2(\cdot), \bar{v}_2(\cdot))\in \mathcal{Q}_2[0,T] \times \mathcal{U}_2[0,T]$  such that
\begin{equation}
\begin{aligned}
&J_2\big(x;\bar{\Theta}_1[\bar{\Theta}_2\bar{x}+\bar{v}_2](\cdot)\bar{x}(\cdot)+\bar{v}_1[\bar{\Theta}_2\bar{x}+\bar{v}_2](\cdot),\bar{\Theta}_2(\cdot)\bar{x}(\cdot)+\bar{v}_2(\cdot)\big)\\
&\leqslant J_2\big(x;\bar{\Theta}_1[\Theta_2\bar{x}^{\Theta_2,v_2}+v_2](\cdot)\bar{x}^{\Theta_2,v_2}(\cdot)+\bar{v}_1[\Theta_2\bar{x}^{\Theta_2,v_2}+v_2](\cdot),\Theta_2(\cdot)\bar{x}^{\Theta_2,v_2}(\cdot)+v_2(\cdot)\big),\\
&\hspace{8cm} \forall\, (\Theta_2(\cdot),v_2(\cdot))\in \mathcal{Q}_2[0,T] \times \mathcal{U}_2[0,T],
\end{aligned}
\end{equation}
where $\bar{x}(\cdot)\equiv\bar{x}^{\bar{\Theta}_2,\bar{v}_2}(\cdot)$ with $\bar{x}^{\Theta_2,v_2}(\cdot)$ being the solution to the {\it closed-loop system} under the {\it closed-loop strategy} $(\Theta_2(\cdot),v_2(\cdot))$ of the leader.

\end{mydef}

\begin{Remark}
We can easily obtain the equation for $\bar{x}^{\Theta_2,v_2}(\cdot)$ in the above definition, by substituting $u_2(\cdot)$ with $\Theta_2(\cdot)\bar{x}^{\Theta_2,v_2}(\cdot)+v_2(\cdot)$ in (\ref{cse}), noting the dependence of $\bar{\Theta}_1(\cdot)$ and $\bar{v}_1(\cdot)$ on $(\Theta_2(\cdot),v_2(\cdot))$. We will give the details in Section 4, when dealing with the problem of the leader. Moreover, we point out that Definition \ref{def2.3} is different from (\ref{problem}), which is one of the main contribution of this paper.
\end{Remark}

\section{LQ problem of the follower}

Let us introduce the following assumptions, which will be in force throughout this paper.
\par \textbf{(H1)} The coefficients of the state equation (\ref{state}) satisfy the following:
\begin{equation}\nonumber
A(\cdot) \in L^1(0,T;\mathbb{R}^{n \times n}),\,\,B_i(\cdot) \in L^2(0,T;\mathbb{R}^{n \times m_i}),\,\,C(\cdot) \in L^2(0,T;\mathbb{R}^{n \times n}),\,\,i=1,2.
\end{equation}
\par \textbf{(H2)} The weighting coefficients in the cost functional (\ref{cfi}) satisfy the following:
\begin{equation}\nonumber
Q_i(\cdot) \in L^1(0,T;\mathbb{S}^n),\,\,R_i(\cdot) \in L^{\infty}(0,T;\mathbb{S}^{m_i})\,\,\mbox{is invertible},\,\,\,G_i \in \mathbb{S}^n,\,\,\,i=1,2.
\end{equation}

\textbf{Problem (SLQ)$_f$}. For any $x \in \mathbb{R}^n$, and given $u_2(\cdot) \in \mathcal{U}_2[0,T]$, find $\bar{u}_1(\cdot)\equiv\bar{u}_1(\cdot;x,u_2) \in \mathcal{U}_1[0,T]$ such that
\begin{equation}\label{follower problem}
J_1(x;\bar{u}_1(\cdot),u_2(\cdot))=\underset{u_1(\cdot)\in \mathcal{U}_1[0,T]} {\min}J_1(x;u_1(\cdot),u_2(\cdot)) \equiv V_1(x).
\end{equation}

It is worth noting that both the open-loop optimal control $\bar{u}_1(\cdot)$ and the value function $V_1(\cdot)$ of the follower depends on the choice of the leader.

First, using the idea of Proposition \ref{opsn}, we are able to obtain the following result.

%follower's open-loop optimal
\begin{mypro}
Let (H1)-(H2) hold. For a given $x \in \mathbb{R}^n$ and $u_2(\cdot) \in \mathcal{U}_2[0,T]$, a state-control pair $(\bar{x}^{u_2}(\cdot),\bar{u}_1(\cdot))$ is an open-loop optimal pair of Problem (SLQ)$_f$ if and only if the following holds:
\begin{equation}\label{fsc}
   B_1(s)^\top \bar{y}(s)+R_1(s)\bar{u}_1(s)=0,\quad a.e.\ s\in[0,T],\ \mathbb{P}\mbox{-}a.s.,
\end{equation}
where $(\bar{y}(\cdot),\bar{z}(\cdot)) \in L^2_{\mathbb{F}}(0,T;\mathbb{R}^n) \times L^2_{\mathbb{F}}(0,T;\mathbb{R}^n)$ is the solution to the following BSDE:
\begin{equation}\left\{
\begin{aligned}
d\bar{y}(s)&=-\big[A(s)^\top \bar{y}(s)+C(s)^\top\bar{z}(s)+Q_1(s)\bar{x}^{u_2}(s)\big]ds+\bar{z}(s)dW(s),\,\quad\,s \in [0,T],\\
 \bar{y}(T)&=G_1\bar{x}^{u_2}(T),
\end{aligned}\right.
\end{equation}
and the following convexity condition holds:
\begin{equation}
\begin{aligned}
&\mathbb{E}\bigg\{ \int_0^T \Big[\big\langle Q_1(s)x_0(s),x_0(s)\big\rangle+\big\langle R_1(s)u_1(s),u_1(s)\big\rangle \Big]dt\\
&\qquad\quad +\big\langle G_1x_0(T),x_0(T)\big\rangle \bigg\} \geqslant 0,\quad \forall u_1(\cdot)\in \mathcal{U}_1[0,T],
\end{aligned}
\end{equation}
where $x_0(\cdot) \in L^2_{\mathbb{F}}(0,T;\mathbb{R}^n)$ is the solution to the following SDE:
\begin{equation}\left\{
\begin{aligned}
dx_0(t)&=\big[A(s)x_0(t)+B_1(s)u_1(s)\big]ds+C(s)x_0(s)dW(s),\quad s \in [0,T],\\
 x_0(0)&=0.
\end{aligned}\right.
\end{equation}
\end{mypro}

Next, take $\Theta_1(\cdot) \in \mathcal{Q}_1[0,T]$ and $v_1(\cdot)\in \mathcal{U}_1[0,T]$. For any $x \in \mathbb{R}^n$ and $u_2(\cdot)\in \mathcal{U}_2[0,T]$, let us consider the closed-loop system (\ref{cse}) and the  corresponding cost functional (\ref{fccf}). The following result characterizes the closed-loop solvability of Problem (SLQ)$_f$.

We will omit some time variables for simplicity, if there is no ambiguity.

%follower closed-loop optimal
\begin{mythm}
Let (H1)-(H2) hold. Then for given $x \in \mathbb{R}^n$ and $u_2(\cdot) \in \mathcal{U}_2[0,T]$, Problem (SLQ)$_f$ admits a closed-loop optimal strategy if and only if the following Riccati equation admits a solution $P^1(\cdot) \in C([0,T];\mathbb{S}^n)$:
\begin{equation}\label{follower Riccati}
\begin{cases}
\dot{P}^1+P^1A+A^\top P^1+C^\top P^1C-P^1B_1R_1^{-1}B_1^\top P^1+Q_1=0,\\
P^1(T)=G_1,\\
R_1 \geqslant 0,\quad a.e.,
\end{cases}
\end{equation}	
and the following BSDE admits a solution $(\eta^{1,u_2}(\cdot),\zeta^{1,u_2}(\cdot)) \in L^2_{\mathbb{F}}(0,T;\mathbb{R}^n) \times L^2_{\mathbb{F}}(0,T;\mathbb{R}^n)$:
\begin{equation}\label{follower BSDE}\left\{
\begin{aligned}
  d\eta^{1,u_2}&=-\big\{ \big[A^\top-P^1B_1R_1^{-1}B_1^\top\big]\eta^{1,u_2}+C^\top\zeta^{1,u_2}+P^1B_2u_2 \big\}ds+\zeta^{1,u_2}dW,\\
\eta^{1,u_2}(T)&=0.
\end{aligned}\right.
\end{equation}	
In this case, the closed-loop optimal strategy $(\bar{\Theta}_1(\cdot),\bar{v}_1(\cdot))$ of Problem (SLQ)$_f$ admits the following representation:
\begin{equation}\label{follower closed-loop optimal}
\begin{cases}
\bar{\Theta}_1=-R_1^{-1} B_1^\top P^1,\\
\bar{v}_1=-R_1^{-1} B_1^\top \eta^{1,u_2}.
\end{cases}
\end{equation}
Further, the value function $V_1(\cdot)$ is given by
\begin{equation}\label{value-1}
\begin{split}
V_1(x)&=\mathbb{E}\biggl\{ \big\langle P^1(0)x,x\big\rangle+2\big\langle \eta^{1,u_2}(0),x\big\rangle+\int_0^T\Big[2\big\langle \eta^{1,u_2},B_2u_2\big\rangle-\big|(R^{-1}_1)^{\frac{1}{2}}B_1^\top\eta^{1,u_2}\big|^2\Big] ds \biggr\}.
\end{split}
\end{equation}
\end{mythm}
The proof here is similar to that in \cite{SLY2016,SY2014}, but for the sake of the integrity of the article, we still give the proof.

{\it Proof.} We first prove the necessity. Given $u_2(\cdot)\in\mathcal{U}_2[0,T]$. Let $(\bar{\Theta}_1(\cdot),\bar{v}_1(\cdot)) \in \mathcal{Q}_1[0,T] \times \mathcal{U}_1[0,T]$ be a closed-loop optimal strategy of Problem (SLQ)$_f$. Then, by Proposition \ref{relation}, $\bar{v}_1(\cdot)$ is an open-loop optimal control of the following LQ problem:
\begin{equation*}\left\{
\begin{aligned}
dx^{\bar{\Theta}_1,v_1,u_2}(s)&=\Big\{\big[A(s)+B_1(s)\bar{\Theta}_1(s)\big]x^{\bar{\Theta}_1,v_1,u_2}(s)+B_1(s)v_1(s)+B_2(s)u_2(s)\Big\}ds\\
                              &\quad+C(s)x^{\bar{\Theta}_1,v_1,u_2}(s)dW(s),\quad s \in [0,T],\\
 x^{\bar{\Theta}_1,v_1,u_2}(0)&=x,
\end{aligned}\right.
\end{equation*}
\begin{equation*}
\begin{aligned}
&\hat{J}_1(x;v_1(\cdot),u_2(\cdot))=\mathbb{E} \bigg\{\int_0^T \Big[ \big\langle \big[Q_1(s)+\bar{\Theta}_1^\top(s) R_1(s)\bar{\Theta}_1(s)\big]x^{\bar{\Theta}_1,v_1,u_2}(s),x^{\bar{\Theta}_1,v_1,u_2}(s) \big\rangle\\
&\ +2\big\langle R_1(s)\bar{\Theta}_1(s)x^{\bar{\Theta}_1,v_1,u_2}(s),v_1(s)\big\rangle+\big\langle R_1(s)v_1(s),v_1(s)\big\rangle \Big]ds
 +\big\langle G_1x^{\bar{\Theta}_1,v_1,u_2}(T),x^{\bar{\Theta}_1,v_1,u_2}(T)\big\rangle \bigg\}.
\end{aligned}
\end{equation*}
Hence, by Proposition \ref{opsn}, for any $x \in \mathbb{R}^n$, the following FBSDE admits a solution triple $(\bar{x}^{u_2}(\cdot),\bar{y}^{u_2}(\cdot),\bar{z}^{u_2}(\cdot)) \in L^2_{\mathbb{F}}(0,T;\mathbb{R}^n)\times L^2_{\mathbb{F}}(0,T;\mathbb{R}^n)\times L^2_{\mathbb{F}}(0,T;\mathbb{R}^n)$:
\begin{equation}\label{follower optimal FBSDE}
\left\{
\begin{aligned}
  d\bar{x}^{u_2}&=\big\{(A+B_1\bar{\Theta}_1)\bar{x}^{u_2}+B_1\bar{v}_1+B_2u_2\big\}ds+C\bar{x}^{u_2}dW,\\
  d\bar{y}^{u_2}&=-\big\{(A+B_1\bar{\Theta}_1)^\top\bar{y}^{u_2}+C^\top \bar{z}^{u_2}+(Q_1+\bar{\Theta}_1^\top R_1\bar{\Theta}_1)\bar{x}^{u_2}+\bar{\Theta}_1^\top R_1\bar{v}_1\big\}ds+\bar{z}^{u_2}dW,\\
\bar{x}^{u_2}(0)&=x,\quad\quad \bar{y}^{u_2}(T)=G_1\bar{x}^{u_2}(T),
\end{aligned}\right.
\end{equation}
with $\bar{x}^{u_2}(\cdot)\equiv x^{\bar{\Theta}_1,\bar{v}_1,u_2}(\cdot)$ and the following stationarity condition holds:
\begin{equation}\label{ follower stationarity condition}
B_1^\top \bar{y}^{u_2}+R_1\bar{\Theta}_1\bar{x}^{u_2}+R_1\bar{v}_1=0, \quad a.e.,\,\,\mathbb{P}\mbox{-}a.s.
\end{equation}
Making use of (\ref{ follower stationarity condition}), we may rewrite the BSDE in (\ref{follower optimal FBSDE}) as follows:
\begin{equation}\nonumber
\begin{aligned}
d\bar{y}^{u_2}&=-\big\{(A+B_1\bar{\Theta}_1)^\top\bar{y}^{u_2}+C^\top \bar{z}^{u_2}+(Q_1+\bar{\Theta}_1^\top R_1\bar{\Theta}_1)\bar{x}^{u_2}+\bar{\Theta}_1^\top R_1\bar{v}_1\big\}ds+\bar{z}^{u_2}dW\\
              &=-\big\{A^\top\bar{y}^{u_2}+C^\top\bar{z}^{u_2}+Q_1\bar{x}^{u_2}+\bar{\Theta}_1^\top(B_1^\top\bar{y}^{u_2}+R_1\bar{\Theta}_1\bar{x}^{u_2}+R_1v_1)\big\}ds+\bar{z}^{u_2}dW\\
              &=-\big\{A^\top\bar{y}^{u_2}+C^\top\bar{z}^{u_2}+Q_1\bar{x}^{u_2}\big\}ds+\bar{z}dW.
\end{aligned}
\end{equation}
Thus, we obtain
\begin{equation}\label{follower simple FBSDE}
\left\{
\begin{aligned}
  d\bar{x}^{u_2}&=\big\{(A+B_1\bar{\Theta}_1)\bar{x}^{u_2}+B_1\bar{v}_1+B_2u_2\big\}ds+C\bar{x}^{u_2}dW,\\
  d\bar{y}^{u_2}&=-\big\{A^\top\bar{y}^{u_2}+C^\top\bar{z}^{u_2}+Q_1\bar{x}^{u_2}\big\}ds+\bar{z}^{u_2}dW,\\
\bar{x}^{u_2}(0)&=x,\quad\quad \bar{y}^{u_2}(T)=G_1\bar{x}^{u_2}(T),\\
                &\hspace{-8mm}B_1^\top \bar{y}^{u_2}+R_1\bar{\Theta}_1\bar{x}^{u_2}+R_1\bar{v}_1=0, \quad a.e.,\,\,\,\mathbb{P}\mbox{-}a.s.
\end{aligned}\right.
\end{equation}
Since the above admits a solution for each $x \in \mathbb{R}^n$, and $(\bar{\Theta}_1(\cdot),\bar{v}_1(\cdot))$ is independent of $x$, by subtraction solutions corresponding $x$ and $0$, the later from the former, we see that for any $x \in \mathbb{R}^n$, the following FBSDE admits an adapted solution $(x(\cdot),y(\cdot),z(\cdot)) \in L^2_{\mathbb{F}}(0,T;\mathbb{R}^n)\times L^2_{\mathbb{F}}(0,T;\mathbb{R}^n)\times L^2_{\mathbb{F}}(0,T;\mathbb{R}^n)$:
\begin{equation}\left\{
\begin{aligned}
dx(s)&=\big\{\big[A(s)+B_1(s)\bar{\Theta}_1(s)\big]x(s)\big\}ds+C(s)x(s)dW(s),\\
dy(s)&=-\big\{A^\top(s)y(s)+C^\top(s)z(s)+Q_1(s)x(s)\big\}ds+z(s)dW,\quad s \in [0,T],\\
 x(0)&=x,\quad\quad y(T)=G_1x(T),\\
     &\hspace{-8mm}B_1^\top(s) y(s)+R_1(s)\bar{\Theta}_1(s)x(s)=0, \quad a.e.\,s \in [0,T],\,\,\mathbb{P}\mbox{-}a.s.
\end{aligned}\right.
\end{equation}
Now, we let
\begin{equation*}\left\{
\begin{aligned}
d\mathbb{X}(s)&=\big\{\big[A(s)+B_1(s)\bar{\Theta}_1(s)\big]\mathbb{X}(s)\big\}ds+C(s)\mathbb{X}(s)dW(s),\quad s \in [0,T],\\
 \mathbb{X}(0)&=I_{n\times n},
\end{aligned}\right.
\end{equation*}
and let
\begin{equation*}\left\{
\begin{aligned}
d\mathbb{Y}(s)&=-\big\{A^\top(s)\mathbb{Y}(s)+C^\top(s)\mathbb{Z}(s)+Q_1(s)\mathbb{X}(s)\big\}ds+\mathbb{Z}(s)dW,\quad s \in [0,T],\\
 \mathbb{Y}(T)&=G_1\mathbb{X}(T).
\end{aligned}\right.
\end{equation*}
Clearly, $\mathbb{X}(\cdot),\,\mathbb{Y}(\cdot),\,\mathbb{Z}(\cdot)$ are all well-defined $\mathbb{S}^n$-matrix valued processes. Further,
\begin{equation}\label{ssc}
B_1^\top(s) \mathbb{Y}(s)+R_1(s)\bar{\Theta}_1(s)\mathbb{X}(s)=0, \quad a.e.\,s \in [0,T],\,\,\mathbb{P}\mbox{-}a.s.
\end{equation}
And $\mathbb{X}(\cdot)^{-1}$ exists, which satisfies the following equation:
\begin{equation}\left\{
\begin{aligned}
d\mathbb{X}^{-1}(s)&=\mathbb{X}^{-1}(s)\big\{C^2(s)-\big[A(s)+B_1(s)\bar{\Theta}_1(s)\big]\big\}ds-\mathbb{X}^{-1}(s)C(s)dW(s),\,\,\,\,s \in [0,T],\\
 \mathbb{X}^{-1}(0)&=I_{n\times n}.
\end{aligned}\right.
\end{equation}
We define
\begin{equation}\nonumber
P^1(\cdot)\triangleq\mathbb{Y}(\cdot)\mathbb{X}(\cdot)^{-1},\qquad \Pi^1(\cdot)\triangleq\mathbb{Z}(\cdot)\mathbb{X}(\cdot)^{-1}.
\end{equation}
Then (\ref{ssc}) implies
\begin{equation}\label{Theta condition}
B_1^\top P^1+R_1\bar{\Theta}_1=0, \quad a.e., \,\,\mathbb{P}\mbox{-}a.s.,
\end{equation}
and thus (since $R_1$ is invertible)
\begin{equation}\label{Theta-1}
\bar{\Theta}_1=-R_1^{-1}B_1^\top P^1, \quad a.e., \,\,\mathbb{P}\mbox{-}a.s.
\end{equation}
Also, by It\^{o}'s formula, we get
\begin{equation}\nonumber
\begin{aligned}
dP^1&=d\mathbb{Y}\mathbb{X}^{-1}=d\mathbb{Y}\cdot\mathbb{X}^{-1}+\mathbb{Y}d\mathbb{X}^{-1}+d\mathbb{Y}\cdot d\mathbb{X}^{-1}\\
    &=\big\{-[A^\top P^1+C^\top \Pi^1+Q_1]+P^1[C^2-A-B_1\bar{\Theta}_1]-\Pi^1C\big\}ds+\big(\Pi^1-P^1C\big)dW.
\end{aligned}
\end{equation}
Let $\Lambda=\Pi^1-P^1C$, then
\begin{equation}\nonumber
\begin{aligned}
dP^1&=\big\{-A^\top P^1-C^\top \Pi^1-Q_1-\Lambda C-P^1A-P^1B_1\bar{\Theta}_1\big\}ds+\Lambda dW\\
    &=-\big\{P^1A+A^\top P^1+\Lambda C+C^\top\Lambda+C^\top P^1C+Q_1+P^1B_1\bar{\Theta}_1\big\}ds+\Lambda dW,
\end{aligned}
\end{equation}
and $P^1(T)=G_1$. Thus, $(P^1(\cdot),\Lambda(\cdot))$ is an adapted solution to a BSDE with deterministic coefficients. Hence, $P^1(\cdot)$ is deterministic and $\Lambda(\cdot)=0$ which means
\begin{equation}
\Pi^1=P^1C.
\end{equation}
Therefore,
\begin{equation}\label{follower RE---}
\dot{P}^1+P^1A+A^\top P^1+C^\top P^1C+P^1B_1\bar{\Theta}_1+Q_1=0.
\end{equation}
Using (\ref{Theta-1}), (\ref{follower RE---}) can be written as
\begin{equation}\label{follower RE}
0=\dot{P}^1+A^\top P^1+P^1A+C^\top P^1C-P^1B_1R_1^{-1}B_1^\top P^1+Q_1.
\end{equation}
Thus, we obtain the Riccati equation in (\ref{follower Riccati}). To determine $\bar{v}_1(\cdot)$, we define
\begin{equation}\nonumber
\eta^{1,u_2}\triangleq \bar{y}^{u_2}-P^1\bar{x}^{u_2},\qquad \zeta^{1,u_2}\triangleq \bar{z}^{u_2}-P^1C\bar{x}^{u_2}.
\end{equation}
Then, noting (\ref{follower simple FBSDE}), (\ref{Theta-1}) and (\ref{follower RE}), we get
\begin{equation}\nonumber
\begin{aligned}
d\eta^{1,u_2}&=d\bar{y}^{u_2}-\dot{P}^1\bar{x}^{u_2}ds-P^1d\bar{x}^{u_2}\\
             &=\big\{-A^\top\bar{y}^{u_2}-C^\top\bar{z}^{u_2}-Q_1\bar{x}^{u_2}+A^\top P^1\bar{x}^{u_2}+P^1A\bar{x}^{u_2}+C^\top P^1C\bar{x}^{u_2}\\
             &\qquad -P^1B_1R_1^{-1}B_1^\top P^1\bar{x}^{u_2}+Q_1\bar{x}^{u_2}-P^1A\bar{x}^{u_2}-P^1B_1\bar{\Theta}_1\bar{x}^{u_2}\\
             &\qquad -P^1B_1\bar{v}_1-P^1B_2u_2\bigr\}ds+(\bar{z}^{u_2}-P^1C\bar{x}^{u_2})dW\\
             &=\big\{-A^\top(\eta^{1,u_2}+P^1\bar{x}^{u_2})-C^\top(\zeta^{1,u_2}+P^1C\bar{x}^{u_2})+A^\top P^1\bar{x}^{u_2}\\
             &\qquad +C^\top P^1C\bar{x}^{u_2}-P^1B_1R_1^{-1}B_1^\top P^1\bar{x}^{u_2}+P^1B_1R_1^{-1}B_1^\top P^1\bar{x}^{u_2}\\
             &\qquad -P^1B_1\bar{v}_1-P^1B_2u_2\bigr\}ds+(\bar{z}^{u_2}-P^1C\bar{x}^{u_2})dW\\
             &=-\big\{ A^\top\eta^{1,u_2}+C^\top\zeta^{1,u_2}+P^1B_1\bar{v}_1+P^1B_2u_2 \big\}ds+\zeta^{1,u_2}dW.
\end{aligned}
\end{equation}
According to (\ref{ follower stationarity condition}), we have
\begin{equation}\label{v_1 condition}
\begin{aligned}
0&=B_1^\top \bar{y}^{u_2}+R_1\bar{\Theta}_1\bar{x}^{u_2}+R_1\bar{v}_1\\
 &=B_1^\top (\eta^{1,u_2}+P^1\bar{x}^{u_2})-R_1R_1^{-1}B_1^\top P^1\bar{x}^{u_2}+R_1\bar{v}_1=B_1^\top\eta^{1,u_2}+R_1\bar{v}_1.
\end{aligned}
\end{equation}
Thus, we have $\bar{v}_1=-R^{-1}_1B_1^\top \eta^{1,u_2}$. Consequently, $(\eta^{1,u_2}(\cdot),\zeta^{1,u_2}(\cdot))$ satisfies the BSDE (\ref{follower BSDE}).

To prove $R_1 \geqslant 0$, as well as the sufficiency, we take any $v_1(\cdot) \in \mathcal{U}_1[0,T]$, and let $x^{\bar{\Theta}_1,v_1,u_2}(\cdot)$ be the corresponding state process. Then, by It\^o's formula, we have
\begin{equation}
\begin{aligned}
&J_1(x;\bar{\Theta}_1(\cdot)x^{\bar{\Theta}_1,v_1,u_2}(\cdot)+v_1(\cdot),u_2(\cdot))\equiv J_1(x;\tilde{u}_1(\cdot),u_2(\cdot))\\
&=\mathbb{E}\biggl\{ \big\langle G_1x^{\bar{\Theta}_1,v_1,u_2}(T),x^{\bar{\Theta}_1,v_1,u_2}(T)\big\rangle
 +\int_0^T\big[\langle Q_1x^{\bar{\Theta}_1,v_1,u_2},x^{\bar{\Theta}_1,v_1,u_2} \rangle +\langle R_1\tilde{u}_1,\tilde{u}_1 \rangle\big] ds \biggr\}\\
&=\mathbb{E}\biggl\{ \big\langle P^1(0)x,x\big\rangle +2\big\langle \eta^{1,u_2}(0),x\big\rangle\\
&\qquad +\int_0^T\Big[\big\langle P^1B_1R^{-1}_1B_1^\top P^1x^{\bar{\Theta}_1,v_1,u_2}+P^1B_1\tilde{u}_1+P^1B_2u_2,x^{\bar{\Theta}_1,v_1,u_2} \big\rangle\\
&\qquad\qquad +\big\langle P^1x^{\bar{\Theta}_1,v_1,u_2},B_1\tilde{u}_1+B_2u_2 \big\rangle +\big\langle R_1\tilde{u}_1,\tilde{u}_1 \big\rangle
 +2\big\langle P^1B_1R^{-1}_1B_1^\top \eta^{1,u_2},x^{\bar{\Theta}_1,v_1,u_2} \big\rangle\\
&\qquad\qquad -2\big\langle P^1B_2u_2,x^{\bar{\Theta}_1,v_1,u_2} \big\rangle +2\big\langle \eta^{1,u_2},B_1\tilde{u}_1+B_2u_2 \big\rangle\Big] ds \biggr\}.
\end{aligned}
\end{equation}
According to (\ref{Theta condition}) and (\ref{v_1 condition}), we have $B^\top_1P^1=-R_1\bar{\Theta}_1$ and $B^\top_1\eta^{1,u_2}=-R_1\bar{v}_1$. Thus,
\begin{equation}
\begin{aligned}
&J_1(x;\bar{\Theta}_1(\cdot)x^{\bar{\Theta}_1,v_1,u_2}(\cdot)+v_1(\cdot),u_2(\cdot))=\mathbb{E}\biggl\{ \big\langle P^1(0)x,x\big\rangle+2\big\langle \eta^{1,u_2}(0),x\big\rangle\\
&\quad +\int_0^T\Big[\big\langle R_1\bar{\Theta}_1x^{\bar{\Theta}_1,v_1,u_2},\bar{\Theta}_1x^{\bar{\Theta}_1,v_1,u_2} \big\rangle-2\big\langle R_1\bar{\Theta}_1x^{\bar{\Theta}_1,v_1,u_2},\tilde{u}_1 \big\rangle
 +\big\langle R_1\tilde{u}_1,\tilde{u}_1 \big\rangle\\
&\qquad\qquad +2\big\langle R_1\bar{\Theta}_1x^{\bar{\Theta}_1,v_1,u_2},\bar{v}_1 \big\rangle+2\big\langle\eta^{1,u_2},B_2u_2 \big\rangle
 -2\big\langle R_1\bar{v}_1,\tilde{u}_1 \big\rangle\Big] ds \biggr\}\\
&=\mathbb{E}\biggl\{ \big\langle P^1(0)x,x\big\rangle+2\big\langle \eta^{1,u_2}(0),x\big\rangle+\int_0^T\Big[2\big\langle \eta^{1,u_2},B_2u_2\big\rangle-\big|(R^{-1}_1)^{\frac{1}{2}}B_1^\top\eta^{1,u_2}\big|^2\Big]ds \biggr\}\\
&\qquad+\mathbb{E}\int_0^T\big\langle R_1(\tilde{u}_1-\bar{\Theta}_1x^{\bar{\Theta}_1,v_1,u_2}-\bar{v}_1),\tilde{u}_1-\bar{\Theta}_1x^{\bar{\Theta}_1,v_1,u_2}-\bar{v}_1\big\rangle ds\\
&=J_1(x;\bar{\Theta}_1(\cdot)\bar{x}^{u_2}(\cdot)+\bar{v}_1(\cdot),u_2(\cdot))+\mathbb{E}\int_0^T\big\langle R_1(v_1-\bar{v}_1),v_1-\bar{v}_1\big\rangle ds.
\end{aligned}
\end{equation}
Hence,
\begin{equation}\nonumber
J_1(x;\bar{\Theta}_1(\cdot)\bar{x}^{u_2}(\cdot)+\bar{v}_1(\cdot),u_2(\cdot)) \leqslant J_1(x;\bar{\Theta}_1(\cdot)x^{\bar{\Theta}_1,v_1,,u_2}(\cdot)+v_1(\cdot),u_2(\cdot)),\quad \forall v_1(\cdot) \in \mathcal{U}_1[0,T],
\end{equation}
if and only if
\begin{equation}\nonumber
R_1 \geqslant 0,\quad a.e.
\end{equation}
In this case, (\ref{value-1}) holds. By Proposition \ref{relation}, it completes the proof. $\qquad\Box$

\section{LQ problem of the leader}

Now, let Problem (SLQ)$_f$ be uniquely closed-loop solvable for given $u_2(\cdot) \in \mathcal{U}_2[0,T]$. Then by (\ref{follower closed-loop optimal}), the follower takes the following closed-loop optimal control:
\begin{equation}\label{follower optimal control}
\begin{aligned}
\bar{u}_1(t)&=\bar{\Theta}_1(t)\bar{x}^{u_2}(t)+\bar{v}_1(t)\\
            &=-R^{-1}_1(t)B^\top_1(t)P^1(t)\bar{x}^{u_2}(t)-R^{-1}_1(t)B^\top_1(t)\eta^{1,u_2}(t),\quad a.e.\,t \in [0,T],
\end{aligned}
\end{equation}
where the process triple $(\bar{x}^{u_2}(\cdot),\eta^{1,u_2}(\cdot),\zeta^{1,u_2}(\cdot))\in L^2_{\mathbb{F}}(0,T;\mathbb{R}^n)\times L^2_{\mathbb{F}}(0,T;\mathbb{R}^n)\times L^2_{\mathbb{F}}(0,T;\mathbb{R}^n)$ satisfies the following FBSDE, which now, is the ``state" equation of the leader:
\begin{equation}\label{leader state}
\left\{
\begin{aligned}
  d\bar{x}^{u_2}&=\big( \hat{A}\bar{x}^{u_2}+\hat{F}_1\eta^{1,u_2}+B_2u_2 \big)ds+C\bar{x}^{u_2}dW\\
         d\eta^{1,u_2}&=-\big( \hat{A}^\top\eta^{1,u_2}+C^\top\zeta^{1,u_2}+\hat{F}_2u_2 \big)ds+\zeta^{1,u_2}dW,\\
\bar{x}^{u_2}(0)&=x,\,\,\,\,\,\eta^{1,u_2}(T)=0.
\end{aligned}
\right.
\end{equation}
In the above, we have denote
$$\hat{A}\triangleq A-B_1R^{-1}_1B_1^\top P^1,\quad \hat{F}_1\triangleq -B_1R^{-1}_1B_1^\top,\quad \hat{F}_2\triangleq P^1B_2.$$

Knowing that the follower has chosen a closed-loop optimal strategy
$$(\bar{\Theta}_1(\cdot),\bar{v}_1(\cdot))\equiv(\bar{\Theta}_1[u_2](\cdot),\bar{v}_1[u_2](\cdot))$$
such that its outcome $\bar{u}_1(\cdot)\equiv\bar{u}_1[u_2](\cdot)$ is of the form (\ref{follower optimal control}), the leader would like to choose an open-loop optimal control $\bar{u}_2(\cdot) \in \mathcal{U}_2[0,T]$ such that the cost functional
\begin{equation}
\begin{aligned}
&\hat{J}_2(x;u_2(\cdot)) \triangleq J_2(x;\bar{u}_1(\cdot),u_2(\cdot))\\
&=\mathbb{E}\bigg\{ \int_0^T \Big[\big\langle Q_2(s)\bar{x}^{u_2}(s),\bar{x}^{u_2}(s)\big\rangle +\big\langle R_2(s)u_2(s),u_2(s)\big\rangle \Big]dt +\big\langle G_2\bar{x}^{u_2}(T),\bar{x}^{u_2}(T)\big\rangle \bigg\}
\end{aligned}
\end{equation}
is minimized. The LQ problem of the leader can be stated as follows.

\vspace{1mm}

\textbf{Problem (SLQ)$_l$}. For given $x \in \mathbb{R}^n$, find a $\bar{u}_2(\cdot) \in \mathcal{U}_2[0,T]$ such that
\begin{equation}\label{LLQ}
\hat{J}_2(x;\bar{u}_2(\cdot))=\underset{u_2(\cdot)\in \mathcal{U}_2[0,T]} {\min}\hat{J}_2(x;u_2(\cdot))\equiv V_2(x).
\end{equation}

Noting that, different from Problem (SLQ)$_f$, the above Problem (SLQ)$_l$ is an LQ stochastic optimal control problem of FBSDE. For its open-loop optimal control $\bar{u}_2(\cdot) \in \mathcal{U}_2[0,T]$, the corresponding process triple $(\bar{x}(\cdot),\bar{\eta}^1(\cdot),\bar{\zeta}^1(\cdot))\equiv (\bar{x}^{\bar{u}_2}(\cdot),\eta^{1,\bar{u}_2}(\cdot),\zeta^{1,\bar{u}_2}(\cdot))$ is called \textit{an open-loop optimal state process triple} and $(\bar{x}(\cdot),\bar{\eta}^1(\cdot),\bar{\zeta}^1(\cdot),\bar{u}_2(\cdot))$ is called \textit{an open-loop optimal quadruple}.

The open-loop solvability of Problem (SLQ)$_l$ can be similarly defined as Definition \ref{def2.1}. And we have the following result first.

%leader open-loop optimal
\begin{mythm}
Let (H1)-(H2) hold. For a given $x \in \mathbb{R}^n$, $(\bar{x}(\cdot),\bar{\eta}^1(\cdot),\bar{\zeta}^1(\cdot),\bar{u}_2(\cdot))$ is an open-loop optimal quadruple of Problem (SLQ)$_l$ if and only if the following stationarity condition holds:
\begin{equation}\label{leader open stationarity}
\hat{F}^\top_2p^{2,\bar{u}_2}+B^\top_2q^{2,\bar{u}_2}+R_2\bar{u}_2=0,\quad a.e.,\, \mathbb{P}\mbox{-}a.s.,
\end{equation}
where $(p^{2,\bar{u}_2}(\cdot),q^{2,\bar{u}_2}(\cdot),k^{2,\bar{u}_2}(\cdot))\in L^2_{\mathbb{F}}(0,T;\mathbb{R}^n)\times L^2_{\mathbb{F}}(0,T;\mathbb{R}^n)\times L^2_{\mathbb{F}}(0,T;\mathbb{R}^n)$ is the solution to the following FBSDE:
\begin{equation}\label{leader adjoint equation}
\left\{
\begin{aligned}
	dp^{2,\bar{u}_2}&=\big(\hat{A}p^{2,\bar{u}_2}+\hat{F}^\top_1 q^{2,\bar{u}_2}\big)ds+Cp^{2,\bar{u}_2}dW,\\
	dq^{2,\bar{u}_2}&=-\big(\hat{A}^\top q^{2,\bar{u}_2}+C^\top k^{2,\bar{u}_2}+Q_2\bar{x}\big)ds+k^{2,\bar{u}_2}dW,\\
  p^{2,\bar{u}_2}(0)&=0,\,\,\,q^{2,\bar{u}_2}(T)=G_2\bar{x}(T),
\end{aligned}
\right.
\end{equation}
and the following convexity condition holds:
\begin{equation}\label{leader convexity condition}
\begin{aligned}
	\mathbb{E} \bigg\{ \big\langle G_2x_{0l}(T),x_{0l}(T)\big\rangle+\int_0^T\Big[ \big\langle Q_2x_{0l},x_{0l} \big\rangle+\big\langle R_2u_2,u_2\big\rangle \Big]ds\bigg\} \geqslant 0,\quad \forall u_2(\cdot) \in \mathcal{U}_2[0,T],
\end{aligned}
\end{equation}
where $(x_{0l}(\cdot),\eta^0(\cdot),\zeta^0(\cdot))\in L^2_{\mathbb{F}}(0,T;\mathbb{R}^n)\times L^2_{\mathbb{F}}(0,T;\mathbb{R}^n)\times L^2_{\mathbb{F}}(0,T;\mathbb{R}^n)$ is the solution to the following FBSDE:
\begin{equation}\label{leader x_0}
\left\{
\begin{aligned}
	dx_{0l}&=\big( \hat{A}x_{0l}+\hat{F}_1\eta^0+B_2u_2 \big)ds+Cx_{0l}dW,\\
	d\eta^0&=-\big( \hat{A}^\top\eta^0+C^\top\zeta^0+\hat{F}_2u_2   \big)ds+\zeta^0dW,\\
  x_{0l}(0)&=0,\,\,\,\,\,\eta^0(T)=0.
\end{aligned}
\right.
\end{equation}
\end{mythm}

\textit{Proof}. Suppose $(\bar{x}^{\bar{u}_2}(\cdot),\eta^{1,\bar{u}_2}(\cdot),\zeta^{1,\bar{u}_2}(\cdot),\bar{u}_2(\cdot))$ is a state-control quadruple corresponding to the given $x \in\mathbb{R}^n$. For any $u_2(\cdot) \in \mathcal{U}_2[0,T]$ and $\epsilon \in \mathbb{R}$, let $u_2^\epsilon(\cdot)=\bar{u}_2(\cdot)+\epsilon u_2(\cdot)$ and $(\bar{x}^{\epsilon}(\cdot)\equiv \bar{x}^{u_2^\epsilon}(\cdot),\eta^{1,\epsilon}(\cdot),\zeta^{1,\epsilon}(\cdot))$ be the corresponding state. Then $(\bar{x}^\epsilon(\cdot),\eta^{1,\epsilon}(\cdot),\zeta^{1,\epsilon}(\cdot))$ satisfies
\begin{equation}\nonumber
\left\{
\begin{aligned}
  d\bar{x}^\epsilon&=\big[ \hat{A}\bar{x}^\epsilon+\hat{F}_1\eta^{1,\epsilon}+B_2(\bar{u}_2+\epsilon u_2) \big]ds+C\bar{x}^\epsilon dW\\
 d\eta^{1,\epsilon}&=-\big[ \hat{A}^\top\eta^{1,\epsilon}+C^\top\zeta^{1,\epsilon}+\hat{F}_2(\bar{u}_2+\epsilon u_2) \big]ds+\zeta^{1,\epsilon} dW,\\
\bar{x}^\epsilon(0)&=x,\,\,\,\,\,\eta^{1,\epsilon}(T)=0.
\end{aligned}
\right.
\end{equation}
Thus, $x_{0l}(\cdot) \equiv \frac{\bar{x}^{\epsilon}(\cdot)-\bar{x}^{\bar{u}_2}(\cdot)}{\epsilon}$ is independent of $\epsilon$ and satisfies (\ref{leader x_0}). Then we get
\begin{equation}\nonumber
\begin{aligned}
&\hat{J}_2(x;\bar{u}_2(\cdot)+\epsilon u_2(\cdot))-\hat{J}_2(x;\bar{u}_2(\cdot))\\
&=2\epsilon\mathbb{E}\biggl\{ \big\langle G_2\bar{x}^{\bar{u}_2}(T),x_{0l}(T)\big\rangle +\int_0^T \Big[ \big\langle  Q_2\bar{x}^{\bar{u}_2},x_{0l} \big\rangle+\big\langle R_2\bar{u}_2,u_2 \big\rangle \Big]ds \biggr\}\\
&\quad +\epsilon^2\mathbb{E}\biggl\{ \big\langle G_2x_{0l}(T),x_{0l}(T)\rangle +\int_0^T \Big[\big\langle Q_2x_{0l},x_{0l} \big\rangle+\big\langle R_2u_2,u_2\big\rangle \Big]ds \biggr\}.
\end{aligned}
\end{equation}
Applying It\^o's formula to $\big\langle q^{2,\bar{u}_2}(\cdot),x_{0l}(\cdot) \big\rangle-\big\langle p^{2,\bar{u}_2}(\cdot),\eta^0(\cdot) \big\rangle$, we obtain
\begin{equation*}
\begin{aligned}
&\hat{J}_2(x;\bar{u}_2(\cdot)+\epsilon u_2(\cdot))-\hat{J}_2(x;\bar{u}_2(\cdot))\\
&=\epsilon\mathbb{E}\biggl\{  \int_0^T \big\langle \hat{F}^\top_2p^{2,\bar{u}_2}+B^\top_2q^{2,\bar{u}_2}+R_2\bar{u}_2,u_2 \big\rangle ds \biggr\}\\
&\quad +\epsilon^2\mathbb{E}\biggl\{ \big\langle G_2x_{0l}(T),x_{0l}(T)\big\rangle+\int_0^T\Big[ \big\langle Q_2x_{0l},x_{0l} \big\rangle+\big\langle R_2u_2,u_2\big\rangle \Big] ds \biggr\}.
\end{aligned}
\end{equation*}
Therefore, $(\bar{x}^{\bar{u}_2}(\cdot),\bar{\eta}^{1,\bar{u}_2}(\cdot),\bar{\zeta}^{1,\bar{u}_2}(\cdot),\bar{u}_2(\cdot))$ is an open-loop optimal quadruple of Problem (SLQ)$_{l}$ if and only if (\ref{leader open stationarity}) and (\ref{leader convexity condition}) hold. The proof is complete. $\qquad\Box$

\vspace{1mm}

Next, as in Definition \ref{def2.3}, we take $\Theta_2(\cdot) \in \mathcal{Q}_2[0,T],\,\,\check{\Theta}_2(\cdot) \in \mathcal{Q}_2[0,T]$ and $v_2(\cdot)\in \mathcal{U}_2[0,T]$. For any $x \in \mathbb{R}^n$, let us consider the following FBSDE:
\begin{equation}\label{leader anticipating closedloop system}
\left\{
\begin{aligned}
  d\bar{x}^{\Theta_2,\check{\Theta}_2,v_2}&=\big[ (\hat{A}+B_2\Theta_2)\bar{x}^{\Theta_2,\check{\Theta}_2,v_2}+(\hat{F}_1
                                           +B_2\check{\Theta}_2)\eta^{1,\Theta_2,\check{\Theta}_2,v_2}+B_2v_2 \big]ds\\
                                          &\quad +C\bar{x}^{\Theta_2,\check{\Theta}_2,v_2}dW,\\
   d\eta^{1,\Theta_2,\check{\Theta}_2,v_2}&=-\big[(\hat{A}+\check{\Theta}_2^\top\hat{F}^\top_2)^\top\eta^{1,\Theta_2,\check{\Theta}_2,v_2}+C^\top\zeta^{1,\Theta_2,\check{\Theta},v_2}\\
                                          &\qquad +\hat{F}_2\Theta_2\bar{x}^{\Theta_2,\check{\Theta}_2,v_2}+\hat{F}_2v_2\big]ds+\zeta^{1,\Theta_2,\check{\Theta}_2,v_2}dW,\\
\bar{x}^{\Theta_2,\check{\Theta}_2,v_2}(0)&=x,\,\,\,\,\,\eta^{1,\Theta_2,\check{\Theta}_2,v_2}(T)=0.
\end{aligned}
\right.
\end{equation}
This is a fully coupled FBSDE which admits a unique solution $(\bar{x}^{\Theta_2,\check{\Theta}_2,v_2},\eta^{1,\Theta_2,\check{\Theta}_2,v_2},\zeta^{1,\Theta_2,\check{\Theta}_2,v_2})\in L^2_{\mathbb{F}}(0,T;\mathbb{R}^n)\times L^2_{\mathbb{F}}(0,T;\mathbb{R}^n)\times L^2_{\mathbb{F}}(0,T;\mathbb{R}^n)$, depending on $\Theta_2(\cdot),\,\check{\Theta}_2(\cdot)$ and $v_2(\cdot)$. (\ref{leader anticipating closedloop system}) is called the {\it closed-loop system} of the original state equation (\ref{leader state}) under the closed-loop strategy $(\Theta_2(\cdot),\check{\Theta}_2(\cdot),\\v_2(\cdot))$ of the leader. Similarly, we point out that $(\Theta_2(\cdot),\check{\Theta}_2(\cdot),v_2(\cdot))$ is independent of the initial state $x$. With the above $\bar{x}^{\Theta_2,\check{\Theta}_2,v_2}(\cdot)$, we define
\begin{equation}\label{laclf}
\begin{aligned}
&\check{J}_2(x;\Theta_2(\cdot)\bar{x}^{\Theta_2,\check{\Theta}_2,v_2}(\cdot)+\check{\Theta}_2(\cdot)\eta^{1,\Theta_2,\check{\Theta}_2,v_2}(\cdot)+v_2(\cdot))\\
&=\mathbb{E} \bigg\{ \big\langle G_2\bar{x}^{\Theta_2,\check{\Theta}_2,v_2}(T),\bar{x}^{\Theta_2,\check{\Theta}_2,v_2}(T)\big\rangle
 +\int_0^T\Big[\big\langle\big[Q_2+\Theta_2^\top R_2\Theta_2\big]\bar{x}^{\Theta_2,\check{\Theta}_2,v_2},\bar{x}^{\Theta_2,\check{\Theta}_2,v_2}\big\rangle\\
&\qquad +\big\langle \check{\Theta}_2^\top R_2 \check{\Theta}_2\eta^{1,\Theta_2,\check{\Theta}_2,v_2},\eta^{1,\Theta_2,\check{\Theta}_2,v_2} \big\rangle+\big\langle R_2v_2,v_2\big\rangle
 +2\big\langle R_2\Theta_2\bar{x}^{\Theta_2,\check{\Theta}_2,v_2},v_2 \big\rangle\\
&\qquad +2\big\langle R_2\Theta_2\bar{x}^{\Theta_2,\check{\Theta}_2,v_2},\check{\Theta}_2\eta^{1,\Theta_2,\check{\Theta}_2,v_2} \big\rangle+2\big\langle R_2\check{\Theta}_2\eta^{1,\Theta_2,\check{\Theta}_2,v_2},v_2\big\rangle \Big]ds\bigg\}.
\end{aligned}
\end{equation}

\begin{mydef}\label{def4.1}
A triple $(\bar{\Theta}_2(\cdot),\bar{\check{\Theta}}_2(\cdot),\bar{v}_2(\cdot)) \in \mathcal{Q}_2[0,T] \times \mathcal{Q}_2[0,T] \times \mathcal{U}_2[0,T]$ is called a \textit{closed-loop optimal strategy} of Problem (SLQ)$_l$ if
\begin{equation}
\begin{aligned}
	&\check{J}_2(x;\bar{\Theta}_2(\cdot)\bar{x}(\cdot)+\bar{\check{\Theta}}_2(\cdot)\bar{\eta}^1(\cdot)+\bar{v}_2(\cdot))\\
    &\leqslant \check{J}_2(x;\Theta_2(\cdot)\bar{x}^{\Theta_2,\check{\Theta}_2,v_2(\cdot)}+\check{\Theta}_2(\cdot)\eta^{1,\Theta_2,\check{\Theta}_2,v_2}(\cdot)+v_2(\cdot)),\\
	&\hspace{1cm}\forall x \in \mathbb{R}^n,\,\,\forall  (\Theta_2(\cdot),\check{\Theta}_2,v_2(\cdot)) \in \mathcal{Q}_2[0,T] \times \mathcal{Q}_2[0,T] \times \mathcal{U}_2[0,T],
\end{aligned}
\end{equation}
where $\bar{x}(\cdot)\equiv \bar{x}^{\bar{\Theta}_2,\bar{\check{\Theta}}_2,\bar{v}_2}(\cdot)$, with $\bar{\eta}^1(\cdot)\equiv \bar{\eta}^{1,\bar{\Theta}_2,\bar{\check{\Theta}}_2,\bar{v}_2}(\cdot)$, $\bar{\zeta}^1(\cdot)\equiv \bar{\zeta}^{1,\bar{\Theta}_2,\bar{\check{\Theta}}_2,\bar{v}_2}(\cdot)$ satisfying (\ref{leader anticipating closedloop system}).
\end{mydef}

The following result is similar to Proposition 3.3 of \cite{SY2014}, and the detailed proof is omitted.
\begin{mypro}\label{pro4.1}
Let (H1)-(H2) hold. Then the following are equivalent:
\par (\romannumeral 1) $(\bar{\Theta}_2(\cdot),\bar{\check{\Theta}}_2(\cdot),\bar{v}_2(\cdot)) \in \mathcal{Q}_2[0,T] \times \mathcal{Q}_2[0,T] \times \mathcal{U}_2[0,T]$ is a closed-loop optimal strategy of Problem (SLQ)$_l$.
\par (\romannumeral 2) The following holds:
\begin{equation*}
\begin{aligned}
   \check{J}_2(x;\bar{\Theta}_2(\cdot)\bar{x}(\cdot)+\bar{\check{\Theta}}_2(\cdot)\bar{\eta}^1(\cdot)+\bar{v}_2(\cdot))
   \leqslant \check{J}_2(x;\bar{\Theta}_2(\cdot)\bar{x}^{\bar{\Theta}_2,\bar{\check{\Theta}}_2,v_2}(\cdot)
  &+\bar{\check{\Theta}}_2(\cdot)\eta^{1,\bar{\Theta}_2,\bar{\check{\Theta}}_2,v_2}(\cdot)+v_2(\cdot)),\\
  &\forall x \in \mathbb{R}^n,\,\, \forall v_2(\cdot) \in \mathcal{U}_2[0,T].
\end{aligned}
\end{equation*}
\par (\romannumeral 3) The following holds:
\begin{equation}\label{open and close equivalent}
\begin{aligned}
	&\check{J}_2(x;\bar{\Theta}_2(\cdot)\bar{x}(\cdot)+\bar{\check{\Theta}}_2(\cdot)\bar\eta^1(\cdot)+\bar{v}_2(\cdot)) \leqslant \check{J}_2(x;u_2(\cdot)),\\
	&\hspace{4cm} \forall x \in \mathbb{R}^n,\,\, \forall u_2(\cdot) \in \mathcal{U}_2[0,T].
\end{aligned}
\end{equation}
\end{mypro}
From (\ref{open and close equivalent}), we can see that for a fixed $x \in \mathbb{R}^n$, the {\it outcome}
\begin{equation}\label{close-eta}
\bar{u}_2(\cdot) \equiv \bar{\Theta}_2(\cdot)\bar{x}(\cdot)+\bar{\check{\Theta}}_2(\cdot)\bar{\eta}^1(\cdot)+\bar{v}_2(\cdot) \in \mathcal{U}_2[0,T]
\end{equation}
of the closed-loop optimal strategy $(\bar{\Theta}_2(\cdot),\bar{\check{\Theta}}_2(\cdot),\bar{v}_2(\cdot))$ is an open-loop optimal control of Problem (SLQ)$_l$. Therefore, Problem (SLQ)$_l$ is closed-loop solvable implies that Problem (SLQ)$_l$ is open-loop solvable.

On the other hand, we can also see that if $(\bar{\Theta}_2(\cdot),\bar{\check{\Theta}}_2(\cdot),\bar{v}_2(\cdot))$ is a closed-loop optimal strategy of Problem (SLQ)$_l$, then $\bar{v}_2(\cdot)$ is an open-loop optimal control of the LQ problem (\ref{leader anticipating closedloop system})-(\ref{laclf}), with $\Theta_2(\cdot)=\bar{\Theta}_2(\cdot),\,\,\check{\Theta}_2(\cdot)=\bar{\check{\Theta}}_2(\cdot)$, which we denote it by {\bf Problem (SLQ)$_{ll}$}. For the open-loop solvability of Problem (SLQ)$_{ll}$, we can similarly obtain the following result.

\begin{mypro}\label{pro4.2}
Let (H1)-(H2) hold. For given $x \in \mathbb{R}^n$, $(\bar{x}(\cdot),\bar{\eta}^1(\cdot),\bar{\zeta}^1(\cdot),\bar{v}_2(\cdot)) \equiv (\bar{x}^{\bar{\Theta}_2,\bar{\check{\Theta}}_2,\bar{v}_2}(\cdot),\\\bar{\eta}^{1,\bar{\Theta}_2,\bar{\check{\Theta}}_2,\bar{v}_2}(\cdot),\bar{\zeta}^{1,\bar{\Theta}_2,\bar{\check{\Theta}}_2,\bar{v}_2}(\cdot),\bar{v}_2(\cdot))$ is an open-loop optimal quadruple of Problem (SLQ)$_{ll}$ if and only if the following stationarity condition holds:
\begin{equation}\label{leader open stationarity---}
\hat{F}^\top_2p^{2,\bar{v}_2}+B^\top_2q^{2,\bar{v}_2}+R_2\bar{\Theta}_2\bar{x}+R_2\bar{\check{\Theta}}_2\bar{\eta}^1+R_2\bar{v}_2=0,\quad a.e.,\, \mathbb{P}\mbox{-}a.s.,
\end{equation}
where $(p^{2,\bar{v}_2}(\cdot),q^{2,\bar{v}_2}(\cdot),k^{2,\bar{v}_2}(\cdot))\in L^2_{\mathbb{F}}(0,T;\mathbb{R}^n)\times L^2_{\mathbb{F}}(0,T;\mathbb{R}^n)\times L^2_{\mathbb{F}}(0,T;\mathbb{R}^n)$ is the solution to the following FBSDE:
\begin{equation}\label{leader adjoint equation---}
\left\{
\begin{aligned}
   dp^{2,\bar{v}_2}&=\bigl[ (\hat{A}+\bar{\check{\Theta}}^\top_2\hat{F}^\top_2)p^{2,\bar{v}_2}+(\hat{F}_1+B_2\bar{\check{\Theta}}_2)^\top q^{2,\bar{v}_2}+\bar{\check{\Theta}}^\top_2R_2\bar{\Theta}\bar{x}\\
                   &\qquad+\bar{\check{\Theta}}^\top_2R_2\bar{\check{\Theta}}_2\bar{\eta}^1+\bar{\check{\Theta}}^\top_2R_2\bar{v}_2 \big]ds+Cp^{2,\bar{v}_2}dW,\\
   dq^{2,\bar{v}_2}&=-\big[ (\hat{A}+B_2\bar{\Theta}_2)^\top q^{2,\bar{v}_2}+C^\top k^{2,\bar{v}_2}+\bar{\Theta}^\top_2\hat{F}^\top_2 p^{2,\bar{v}_2} +\bar{\Theta}^\top_2R_2\bar{\check{\Theta}}_2\bar{\eta}^1\\
                   &\qquad\ +(Q_2+\bar{\Theta}^\top_2R_2\bar{\Theta}_2)\bar{x}+\bar{\Theta}^\top_2R_2\bar{v}_2 \big]ds+k^{2,\bar{v}_2}dW,\\
 p^{2,\bar{v}_2}(0)&=0,\,\,\,q^{2,\bar{v}_2}(T)=G_2\bar{x}(T),
\end{aligned}
\right.
\end{equation}
and the following convexity condition holds:
\begin{equation}\label{leader convexity condition---}
\begin{aligned}
&\mathbb{E} \bigg\{ \big\langle G_2x^{v_2}_{0l}(T),x^{v_2}_{0l}(T)\big\rangle
+\int_0^T\Big[\big\langle\big[Q_2+\bar{\Theta}_2^\top R_2\bar{\Theta}_2\big]x^{v_2}_{0l},x^{v_2}_{0l}\big\rangle+2\big\langle R_2\bar{\Theta}_2x^{v_2}_{0l},\bar{\check{\Theta}}_2\eta^{0,v_2} \big\rangle\\
&\quad+2\big\langle R_2\bar{\Theta}_2x^{v_2}_{0l},v_2 \big\rangle+\big\langle \bar{\check{\Theta}}_2^\top R_2 \bar{\check{\Theta}}_2\eta^{0,v_2},\eta^{0,v_2}\big\rangle+2\big\langle R_2\bar{\check{\Theta}}_2\eta^{0,v_2},v_2\big\rangle+\big\langle R_2v_2,v_2\big\rangle \Big]ds\bigg\} \geqslant 0,\\
&\hspace{10cm}\quad \forall v_2(\cdot)\in\mathcal{U}_2[0,T],
\end{aligned}
\end{equation}
where $(x_{0l}^{v_2}(\cdot),\eta^{0,v_2}(\cdot),\zeta^{0,v_2}(\cdot))\in L^2_{\mathbb{F}}(0,T;\mathbb{R}^n)\times L^2_{\mathbb{F}}(0,T;\mathbb{R}^n)\times L^2_{\mathbb{F}}(0,T;\mathbb{R}^n)$ is the solution to the following FBSDE:
\begin{equation}\label{leader x_0---}
\left\{
\begin{aligned}
   dx_{0l}^{v_2}&=\big[ (\hat{A}+B_2\bar{\Theta}_2)x_{0l}^{v_2}+(\hat{F}_1+B_2\bar{\check{\Theta}}_2)\eta^{0,v_2}+B_2v_2 \big]ds+Cx_{0l}^{v_2}dW,\\
   d\eta^{0,v_2}&=-\big[(\hat{A}^\top+\hat{F}_2\bar{\check{\Theta}}_2)\eta^{0,v_2}+C^\top\zeta^{0,v_2}+\hat{F}_2\bar{\Theta}_2x_{0l}^{v_2}+\hat{F}_2v_2\big]ds+\zeta^{0,v_2}dW,\\
 x_{0l}^{v_2}(0)&=0,\,\,\,\,\,\eta^{0,v_2}(T)=0.
\end{aligned}
\right.
\end{equation}
\end{mypro}

\begin{Remark}
The introduction of Problem (SLQ)$_{ll}$ above, is to give the characterization of the closed-loop optimal strategy of the leader for Problem (SLQ)$_l$. However, if we consider the outcome of the closed-loop strategy of the form $u_2(\cdot)=\Theta_2(\cdot)\bar{x}^{\Theta_2,\check{\Theta}_2,v_2}(\cdot)+\check{\Theta}_2(\cdot)\eta^{1,\Theta_2,\check{\Theta}_2,v_2}(\cdot)+v_2(\cdot)$ as (\ref{close-eta}), it is anticipating since $\eta^{1,\Theta_2,\check{\Theta}_2,v_2}(\cdot)$ exists. This is not realistic. We point out that we overcome this difficulty inspired by Yong \cite{Yong2002}, to give some necessary conditions for the nonanticipating closed-loop optimal strategy of the leader for Problem (SLQ)$_l$. This is one of the main contributions of this paper.
\end{Remark}

Instead of (\ref{leader anticipating closedloop system}), we consider the following closed-loop system:
\begin{equation}\label{leader nonanticipating closed-loop state}
\left\{
\begin{aligned}
  d\bar{x}^{\Theta_2,\tilde{\Theta}_2,v_2}&=\big[ (\hat{A}+B_2\Theta_2)\bar{x}^{\Theta_2,\tilde{\Theta},v_2}+\hat{F}_1\eta^{1,\Theta_2,\tilde{\Theta},v_2}+B_2\tilde{\Theta}_2p^2+B_2v_2 \big]ds\\
                                          &\quad +C\bar{x}^{\Theta_2,\tilde{\Theta}_2,v_2}dW,\\
     d\eta^{1,\Theta_2,\tilde{\Theta},v_2}&=-\big[\hat{A}^\top\eta^{1,\Theta_2,\tilde{\Theta},v_2}+C^\top\zeta^{1,\Theta_2,\tilde{\Theta},v_2}+\hat{F}_2\Theta_2\bar{x}^{\Theta_2,\tilde{\Theta}_2,v_2}\\
                                          &\qquad +\hat{F}_2\tilde{\Theta}_2p^2+\hat{F}_2v_2\big]ds+\zeta^{1,\Theta_2,\tilde{\Theta},v_2}dW,\\
\bar{x}^{\Theta_2,\tilde{\Theta}_2,v_2}(0)&=x,\,\,\,\,\,\eta^{1,\Theta_2,\tilde{\Theta},v_2}(T)=0,
\end{aligned}
\right.
\end{equation}	
with the cost functional
\begin{equation}\label{lnclf}
\begin{aligned}
&\hat{J}_2(x;\Theta_2(\cdot)\bar{x}^{\Theta_2,\tilde{\Theta}_2,v_2}(\cdot)+\tilde{\Theta}_2(\cdot)p^2(\cdot)+v_2(\cdot))=\mathbb{E} \bigg\{ \big\langle G_2\bar{x}^{\Theta_2,\tilde{\Theta}_2,v_2}(T),\bar{x}^{\Theta_2,\tilde{\Theta}_2,v_2}(T)\big\rangle\\
&\quad +\int_0^T\Big[\big\langle\big[Q_2+\Theta_2^\top R_2\Theta_2\big]\bar{x}^{\Theta_2,\tilde{\Theta}_2,v_2},\bar{x}^{\Theta_2,\tilde{\Theta}_2,v_2}\big\rangle+\big\langle \check{\Theta}_2^\top R_2 \tilde{\Theta}_2p^2,p^2 \big\rangle+\big\langle R_2v_2,v_2\big\rangle\\
&\qquad \quad +2\big\langle R_2\Theta_2\bar{x}^{\Theta_2,\tilde{\Theta}_2,v_2},v_2 \big\rangle+2\big\langle R_2\Theta_2\bar{x}^{\Theta_2,\tilde{\Theta}_2,v_2},\tilde{\Theta}_2p^2 \big\rangle+2\big\langle R_2\tilde{\Theta}_2p^2,v_2\big\rangle \Big]ds\bigg\}.
\end{aligned}
\end{equation}	
where $p^2(\cdot)\in L^2_\mathbb{F}(0,T;\mathbb{R}^n)$ with $(q^2(\cdot),k^2(\cdot))\in L^2_\mathbb{F}(0,T;\mathbb{R}^n)\times  L^2_\mathbb{F}(0,T;\mathbb{R}^n)$ is the solution to the following adjoint FBSDE:
\begin{equation}\label{leader nonanticipating adjoint functional}
\left\{
\begin{aligned}
  dp^2&=\bigl( \hat{A}p^2+\hat{F}^\top_1 q^2 \big)ds+Cp^2dW,\\
  dq^2&=-\big[ (\hat{A}+B_2\bar{\Theta}_2)^\top q^2+C^\top k^2+(\bar{\Theta}^\top_2\hat{F}^\top_2+\bar{\Theta}^\top_2R_2\bar{\tilde{\Theta}}_2)p^2\\
      &\qquad+(Q_2+\bar{\Theta}^\top_2R_2\bar{\Theta}_2)\bar{x}+\bar{\Theta}^\top_2R_2\bar{v}_2 \big]ds+k^2dW,\\
p^2(0)&=0,\quad q^2(T)=G_2\bar{x}(T).
\end{aligned}
\right.
\end{equation}
Moreover, the following stationary condition holds
\begin{equation}\label{leader nonanticipating stationary condition}
(\hat{F}^\top_2+R_2\bar{\tilde{\Theta}}_2) p^2+B^\top_2q^2+R_2\bar{\Theta}_2\bar{x}+R_2\bar{v}_2=0,\quad a.e.,\, \mathbb{P}\mbox{-}a.s.,
\end{equation}	
where $(\bar{x}(\cdot),\bar{\eta}^1(\cdot),\bar{\zeta}^1(\cdot)) \equiv (\bar{x}^{\bar{\Theta}_2,\bar{\tilde{\Theta}}_2,\bar{v}_2}(\cdot),\eta^{1,\bar{\Theta}_2,\bar{\tilde{\Theta}}_2,\bar{v}_2}(\cdot),\zeta^{1,\bar{\Theta}_2,\bar{\tilde{\Theta}}_2,\bar{v}_2}(\cdot))$ is the optimal quadruple of the problem (\ref{leader nonanticipating closed-loop state})-(\ref{lnclf}). Similarly, we can give the equivalent definitions of the closed-loop optimal strategy as Definition \ref{def4.1} and Proposition \ref{pro4.1}.

Making use of the stationary condition in (\ref{leader nonanticipating stationary condition}), we may rewrite the BSDE in (\ref{leader nonanticipating adjoint functional}) as follows:
\begin{equation}\label{simple BSDE}
\begin{aligned}
dq^2&=-\big[ (\hat{A}+B_2\bar{\Theta}_2)^\top q^2+C^\top k^2+(\bar{\Theta}^\top_2\hat{F}^\top_2+\bar{\Theta}^\top_2R_2\bar{\tilde{\Theta}}_2)p^2+(Q_2+\bar{\Theta}^\top_2R_2\bar{\Theta}_2)\bar{x}\\
    &\qquad +\bar{\Theta}^\top_2R_2\bar{v}_2 \big]ds+k^{2,\bar{v}_2}dW,\\
    &=-\big[ \hat{A}^\top q^2+C^\top k^2+Q_2\bar{x}+\bar{\Theta}_2^\top (B^\top_2q^2+\hat{F}^\top_2 p^2+R_2\bar{\tilde{\Theta}}_2p^2+R_2\bar{\Theta}_2\bar{x}\\
    &\qquad +R_2\bar{v}_2) \big]ds+k^2dW\\
    &=-\big( \hat{A}^\top q^2+C^\top k^2+Q_2\bar{x}\big)ds+k^2dW.
\end{aligned}
\end{equation}
For convenience, we write the state equations (\ref{leader nonanticipating closed-loop state}) and the adjoint equations (\ref{leader nonanticipating adjoint functional}) together (noting (\ref{simple BSDE})), to obtain
\begin{equation}\label{leader closedloop optimal system simple}
\left\{
\begin{aligned}
d\bar{x}&=\big[ (\hat{A}+B_2\bar{\Theta}_2)\bar{x}+\hat{F}_1\bar{\eta}^1+B_2\bar{\tilde{\Theta}}_2p^2+B_2\bar{v}_2 \big]ds+C\bar{x}dW,\\
d\bar{\eta}^1&=-\big( \hat{A}^\top\bar{\eta}^1+C^\top\bar{\zeta}^1+\hat{F}_2\bar{\Theta}_2\bar{x}+\hat{F}_2\bar{\tilde{\Theta}}_2p^2+\hat{F}_2\bar{v}_2 \big)ds+\bar{\zeta}^1dW,\\
dp^2&=\bigl( \hat{A}p^2+\hat{F}^\top_1 q^2 \big)ds+Cp^2dW,\\
dq^2&=-\big( \hat{A}^\top q^2+C^\top k^2+Q_2\bar{x}\big)ds+k^2dW,\\
\bar{x}(0)&=x,\,\,\,\eta^1(T)=0,\,\,\,p^2(0)=0,\,\,\,q^2(T)=G_2\bar{x}(T),\\
&\hspace{-8mm}(\hat{F}^\top_2+R_2\bar{\tilde{\Theta}}_2) p^2+B^\top_2q^2+R_2\bar{\Theta}_2\bar{x}+R_2\bar{v}_2=0,\quad a.e.,\, \mathbb{P}\mbox{-}a.s.,
\end{aligned}
\right.
\end{equation}
Note that the above is a coupled FBSDEs system which is further coupled through the last relation. Next, let us set
\begin{equation}
X \triangleq\left(\begin{matrix} \bar{x}  \\ p^2  \end{matrix}\right),\,\,\,Y \triangleq\left(\begin{matrix} q^2  \\ \bar{\eta}^1  \end{matrix}\right),\,\,\,
Z \triangleq\left(\begin{matrix} k^2  \\ \bar{\zeta}^1  \end{matrix}\right),\,\,\,\boldsymbol{\Theta_2}\triangleq\left(\begin{matrix} \bar{\Theta}_2 & \bar{\tilde{\Theta}}_2 \end{matrix}\right),
\end{equation}
and
\begin{equation*}
\begin{cases}
\mathcal{A} \triangleq\left(\begin{matrix} \hat{A} & 0 \\ 0 & \hat{A} \end{matrix}\right),\quad
\mathcal{B}_2 \triangleq\left(\begin{matrix} B_2 \\ 0 \end{matrix}\right),\quad
\mathcal{F}_1 \triangleq\left(\begin{matrix} 0 & \hat{F}_1 \\ \hat{F}_1^\top&0 \end{matrix}\right),\quad
\mathcal{C} \triangleq\left(\begin{matrix} C & 0  \\ 0 & C \end{matrix}\right),\\
\mathcal{Q}_2 \triangleq\left(\begin{matrix} Q_2 & 0 \\ 0 & 0 \end{matrix}\right),\quad
\mathcal{F}_2 \triangleq\left(\begin{matrix} 0 \\ \hat{F}_2 \end{matrix}\right),\quad
\mathcal{G}_2 \triangleq\left(\begin{matrix} G_2&0 \\ 0&0 \end{matrix}\right),\quad
X_0 \triangleq\left(\begin{matrix} x  \\ 0  \end{matrix}\right).
\end{cases}
\end{equation*}
Then (\ref{leader closedloop optimal system simple}) is equivalent to the following FBSDE:
\begin{equation}\label{High dimension optimal system simple}
\begin{cases}
dX=\big[ (\mathcal{A}+\mathcal{B}_2\boldsymbol{\Theta}_2)X+\mathcal{F}_1Y+\mathcal{B}_2\bar{v}_2  \big]ds+\mathcal{C}XdW,\\
dY=-\big[ (\mathcal{Q}_2+\mathcal{F}_2\boldsymbol{\Theta}_2)X+\mathcal{A}^\top Y+\mathcal{C}^\top Z+\mathcal{F}_2\bar{v}_2 \big]ds+ZdW,\\
X(0)=X_0,\,\,\,Y(T)=\mathcal{G}_2X(T),
\end{cases}
\end{equation}
whose solution triple $(X(\cdot),Y(\cdot),Z(\cdot))\in L^2_{\mathbb{F}}(0,T;\mathbb{R}^{2n})\times L^2_{\mathbb{F}}(0,T;\mathbb{R}^{2n})\times L^2_{\mathbb{F}}(0,T;\mathbb{R}^{2n})$, together with the following condition holds:
\begin{equation}\label{High dimension stationary}
(R_2\boldsymbol{\Theta}_2+\mathcal{F}^\top_2)X+\mathcal{B}^\top_2Y+R_2\bar{v}_2=0,\quad a.e.,\, \mathbb{P}\mbox{-}a.s.
\end{equation}

For the closed-loop optimal strategies of the leader, we have the following result.
\begin{mythm}\label{thm4.2}
Let (H1)-(H2) hold, if Problem (SLQ)$_l$ is closed-loop solvable, then the closed-loop optimal strategy $(\boldsymbol\Theta_2(\cdot),\bar{v}_2(\cdot))\equiv(\bar{\Theta}_2(\cdot),\bar{\tilde{\Theta}}_2(\cdot),\bar{v}_2(\cdot))\in \mathcal{Q}_2[0,T] \times \mathcal{Q}_2[0,T] \times \mathcal{U}_2[0,T]$ admits the following representation:
\begin{equation}\label{close-loop leader}
\left\{
\begin{aligned}
        \bar{\Theta}_2&=-R^{-1}_2B^\top_2P_1,\\
\bar{\tilde{\Theta}}_2&=-R^{-1}_2(B^\top_2P_2+\hat{F}^\top_2),\\
             \bar{v}_2&=0,\qquad\qquad a.e.,\, \mathbb{P}\mbox{-}a.s.
\end{aligned}
\right.
\end{equation}
where $P(\cdot)\equiv\left(\begin{matrix} P_1(\cdot)&P_2(\cdot) \\ P_2(\cdot)^\top&P_4(\cdot) \end{matrix}\right) \in C([0,T];\mathbb{S}^{2n \times 2n})$ is the solution to the following Riccati equation:
\begin{equation}\label{leader RE}
\begin{cases}
&\dot{P}+\mathcal{A}^\top P+P\mathcal{A}+\mathcal{C}^\top P\mathcal{C}+P\mathcal{F}_1P+\mathcal{Q}_2-(P\mathcal{B}_2+\mathcal{F}_2)R^{-1}_2(\mathcal{B}^\top_2P+\mathcal{F}^\top_2)=0,\\
&P(T)=\mathcal{G}_2.
\end{cases}
\end{equation}
In this case, the closed-loop optimal control of the leader is $\bar{u}_2(\cdot)=\boldsymbol\Theta_2(\cdot)X(\cdot)$, where $X(\cdot)\in L^2_{\mathbb{F}}(0,T;\mathbb{R}^{2n})$ is the solution to the following SDE:
\begin{equation}\label{nonanticipating X}\left\{
\begin{aligned}
  dX&=\big[\mathcal{A}-\mathcal{B}_2R^{-1}_2(\mathcal{B}^\top_2P+\mathcal{F}^\top_2)+\mathcal{F}_1P\big]Xds+\mathcal{C}XdW,\\
X(0)&=X_0.
\end{aligned}\right.
\end{equation}
Further, the value function of the leader admits the following representation:
\begin{equation}\label{value function-leader}
V_2(x)=\big\langle P_1(0)x,x \big\rangle.
\end{equation}
\end{mythm}

\textit{Proof.} Let $(\bar{\Theta}_2(\cdot),\bar{\tilde{\Theta}}_2(\cdot),\bar{v}_2(\cdot))$ be a closed-loop optimal strategy of Probelm (SLQ)$_l$. Since (\ref{High dimension optimal system simple}) admits a solution for each $X_0 \in \mathbb{R}^{2n}$, and $(\boldsymbol{\Theta}_2(\cdot),\bar{v}_2(\cdot))$ is independent of $x$, by substracting solutions corresponding $X_0$ and $0$, the later from the former, we see that for any $X_0 \in \mathbb{R}^{2n}$, the following FBSDE admits an adapted solution $(\tilde{X}(\cdot),\tilde{Y}(\cdot),\tilde{Z}(\cdot))$:
\begin{equation*}
\left\{
\begin{aligned}
  d\tilde{X}&=\big[ (\mathcal{A}+\mathcal{B}_2\boldsymbol{\Theta}_2)\tilde{X}+\mathcal{F}_1\tilde{Y} \big]ds+\mathcal{C}\tilde{X}dW,\\
  d\tilde{Y}&=-\big[ (\mathcal{Q}_2+\mathcal{F}_2\boldsymbol{\Theta}_2)\tilde{X}+\mathcal{A}^\top \tilde{Y}+\mathcal{C}^\top \tilde{Z} \big]ds+\tilde{Z}dW,\\
\tilde{X}(0)&=X_0,\,\,\,\tilde{Y}(T)=\mathcal{G}_2\tilde{X}(T).
\end{aligned}
\right.
\end{equation*}
Now, we let
\begin{equation}
\left\{
\begin{aligned}
  d\mathbb{X}&=\big[ (\mathcal{A}+\mathcal{B}_2\boldsymbol{\Theta}_2)\mathbb{X}+\mathcal{F}_1\mathbb{Y} \big]ds+\mathcal{C}\mathbb{X}dW,\\
  d\mathbb{Y}&=-\big[ (\mathcal{Q}_2+\mathcal{F}_2\boldsymbol{\Theta}_2)\mathbb{X}+\mathcal{A}^\top \mathbb{Y}+\mathcal{C}^\top \mathbb{Z}\big]ds+\mathbb{Z}dW,\\
\mathbb{X}(0)&=I_{2n\times 2n},\,\,\,\mathbb{Y}(T)=\mathcal{G}_2\mathbb{X}(T).
\end{aligned}
\right.
\end{equation}
Clearly, $\mathbb{X}(\cdot),\mathbb{Y}(\cdot),\mathbb{Z}(\cdot)$ are all well-defined $(2n\times 2n)$-matrix valued processes. Further, (\ref{High dimension stationary}) is equivalent to
\begin{equation}\label{simply stationary}
(R_2\boldsymbol{\Theta}_2+\mathcal{F}^\top_2)\mathbb{X}+\mathcal{B}^\top_2\mathbb{Y}=0,\quad a.e.,\, \mathbb{P}\mbox{-}a.s.
\end{equation}

Drawing on the method of Yong \cite{Yong2006}, we can check that $\mathbb{X}(\cdot)^{-1}$ exists and satisfies the following SDE:
\begin{equation}
\left\{
\begin{aligned}
  d\mathbb{X}^{-1}&=\Big\{-\mathbb{X}^{-1}\big[(\mathcal{A}+\mathcal{B}_2\boldsymbol{\Theta}_2)\mathbb{X}+\mathcal{F}_1\mathbb{Y}\big]\mathbb{X}^{-1}+\mathbb{X}^{-1}\mathcal{C}^2 \Big\}ds-\mathbb{X}^{-1}\mathcal{C}dW,\\
\mathbb{X}(0)^{-1}&=I_{2n\times 2n}.
\end{aligned}
\right.
\end{equation}
We define
\begin{equation}\label{PX-1}
P(\cdot)\triangleq\mathbb{Y}(\cdot)\mathbb{X}(\cdot)^{-1},\qquad \Pi(\cdot)\triangleq\mathbb{Z}(\cdot)\mathbb{X}(\cdot)^{-1}.
\end{equation}
By It\^o's formula, we obtain
\begin{equation}\nonumber
\begin{aligned}
dP&=\Big\{ -\big[(\mathcal{Q}_2+\mathcal{F}_2\boldsymbol{\Theta}_2)\mathbb{X}+\mathcal{A}^\top \mathbb{Y}+\mathcal{C}^\top \mathbb{Z}\big]\mathbb{X}^{-1}
   -\mathbb{Y}\mathbb{X}^{-1}(\mathcal{A}+\mathcal{B}_2\boldsymbol{\Theta}_2)\mathbb{X}\mathbb{X}^{-1}\\
  &\qquad -\mathbb{Y}\mathbb{X}^{-1}\mathcal{F}_1\mathbb{Y}\mathbb{X}^{-1}+\mathbb{Y}\mathbb{X}^{-1}\mathcal{C}^2-\mathbb{Z}\mathbb{X}^{-1}\mathcal{C}\Big\}ds
   +\big( \mathbb{Z}\mathbb{X}^{-1}-\mathbb{Y}\mathbb{X}^{-1}\mathcal{C} \big)dW\\
  &=\big[ -(\mathcal{Q}_2+\mathcal{F}_2\boldsymbol{\Theta}_2)-\mathcal{A}^\top P-\mathcal{C}^\top \Pi-P(\mathcal{A}+\mathcal{B}_2\boldsymbol{\Theta}_2)\\
  &\qquad -P\mathcal{F}_1P+P\mathcal{C}^2 -\Pi\mathcal{C} \big]ds+(\Pi-P\mathcal{C})dW.
\end{aligned}
\end{equation}
Let
\begin{equation}\nonumber
\Lambda(\cdot)\triangleq\Pi(\cdot)-P(\cdot)\mathcal{C}(\cdot)
\end{equation}
which leads to
\begin{equation}\nonumber
\begin{aligned}
dP&=\big[ -(\mathcal{Q}_2+\mathcal{F}_2\boldsymbol{\Theta}_2)-\mathcal{A}^\top P-\mathcal{C}^\top (\Lambda+P\mathcal{C})-P(\mathcal{A}+\mathcal{B}_2\boldsymbol{\Theta}_2) \\
  &\qquad -P\mathcal{F}_1P+P\mathcal{C}^2 -(\Lambda+P\mathcal{C})\mathcal{C} \big]ds +\Lambda dW\\
  &=-\big[ (\mathcal{Q}_2+\mathcal{F}_2\boldsymbol{\Theta}_2)+\mathcal{A}^\top P+\mathcal{C}^\top (\Lambda+P\mathcal{C})+P(\mathcal{A}+\mathcal{B}_2\boldsymbol{\Theta}_2)+P\mathcal{F}_1P+\Lambda\mathcal{C} \big]ds+\Lambda dW,
\end{aligned}
\end{equation}
and $P(T)=\mathcal{G}_2$. Thus, $(P(\cdot),\Lambda(\cdot))$ is the adapted solution to a BSDE with deterministic coefficients. Hence, $P(\cdot)$ is deterministic and $\Lambda(\cdot)=0$ which means
\begin{equation}\nonumber
\Pi(\cdot)=P(\cdot)\mathcal{C}(\cdot).
\end{equation}
Therefore, we get
\begin{equation}\label{leader riccati}
\dot{P}+\mathcal{A}^\top P+P\mathcal{A}+\mathcal{C}^\top P\mathcal{C}+P\mathcal{F}_1P+\mathcal{Q}_2+(P\mathcal{B}_2+\mathcal{F}_2)\boldsymbol{\Theta}_2=0.
\end{equation}
Moreover, (\ref{simply stationary}) and (\ref{PX-1}) imply
\begin{equation}\label{leader Theta condition}
R_2\boldsymbol{\Theta}_2+\mathcal{F}^\top_2+\mathcal{B}^\top_2P=0,\quad a.e.,\, \mathbb{P}\mbox{-}a.s.
\end{equation}
Thus
\begin{equation}\label{Theta2}
\boldsymbol{\Theta}_2=-R^{-1}_2(\mathcal{B}^\top_2P+\mathcal{F}^\top_2),
\end{equation}
and since
\begin{equation*}
\begin{aligned}
\boldsymbol{\Theta}_2& \equiv\left(\begin{matrix} \bar{\Theta}_2 & \bar{\tilde{\Theta}}_2\end{matrix}\right)
                =-R^{-1}_2\bigg[ \left(\begin{matrix} B^\top_2 & 0\end{matrix}\right)\left(\begin{matrix} P_1 & P_2\\P^\top_2&P_4\end{matrix}\right)
                 +\left(\begin{matrix} 0 & \hat{F}^\top_2\end{matrix}\right) \bigg]\\
              & =\left(\begin{matrix} -R^{-1}_2B^\top_2P_1 & -R^{-1}_2(B^\top_2P_2+\hat{F}^\top_2)\end{matrix}\right),
\end{aligned}
\end{equation*}
we have
\begin{equation}\label{theta2}
\bar{\Theta}_2=-R^{-1}_2B^\top_2P_1,\qquad \bar{\tilde{\Theta}}_2=-R^{-1}_2(B^\top_2P_2+\hat{F}^\top_2).
\end{equation}

Plugging the above into (\ref{leader riccati}), we obtain the Riccati equation in (\ref{leader RE}). To determine $\bar{v}_2(\cdot)$, we define
\begin{equation}\nonumber
\eta(\cdot)\triangleq Y(\cdot)-P(\cdot)X(\cdot),\qquad \zeta(\cdot)\triangleq Z(\cdot)-P(\cdot)\mathcal{C}(\cdot)X(\cdot).
\end{equation}
Consequently,
\begin{equation}\label{eta}
\begin{aligned}
d\eta&=\big[ -(\mathcal{Q}_2+\mathcal{F}_2\boldsymbol{\Theta}_2)X-\mathcal{A}^\top Y-\mathcal{C}^\top Z-\mathcal{F}_2\bar{v}_2+   \mathcal{A}^\top PX\\
     &\qquad +P\mathcal{A}X+\mathcal{C}^\top P\mathcal{C}X+P\mathcal{B}_2\boldsymbol{\Theta}_2X+P\mathcal{F}_1PX+(\mathcal{Q}_2+\mathcal{F}_2\boldsymbol{\Theta}_2)X\\
     &\qquad -P(\mathcal{A}+\mathcal{B}_2\boldsymbol{\Theta}_2)X-P\mathcal{F}_1Y-P  \mathcal{B}_2\bar{v}_2 \big]ds+\big(Z-P\mathcal{C}X\big)dW\\
     &=-\big[ \mathcal{A}^\top (\eta+PX)+\mathcal{C}^\top (\zeta+P\mathcal{C}X)+\mathcal{F}_2\bar{v}_2 -\mathcal{A}^\top PX-P\mathcal{A}X-\mathcal{C}^\top P\mathcal{C}X\\
     &\qquad -P\mathcal{B}_2\boldsymbol{\Theta}_2X-P\mathcal{F}_1PX+P(\mathcal{A}+\mathcal{B}_2\boldsymbol{\Theta}_2)X+P\mathcal{F}_1(\eta+PX)+P\mathcal{B}_2\bar{v}_2 \big]ds+\zeta dW\\
     &=-\big[ (\mathcal{A}^\top+P\mathcal{F}_1)\eta+\mathcal{C}^\top\zeta+(\mathcal{F}_2+P\mathcal{B}_2)\bar{v}_2 \big]ds+\zeta dW.
\end{aligned}
\end{equation}
According to (\ref{High dimension stationary}) and (\ref{leader Theta condition}), we have
\begin{equation}\nonumber
\begin{aligned}
0&=(R_2\boldsymbol{\Theta}_2+\mathcal{F}^\top_2)X+\mathcal{B}^\top_2(\eta+PX)+R_2\bar{v}_2\\
 &=(R_2\boldsymbol{\Theta}_2+\mathcal{F}^\top_2+\mathcal{B}^\top_2P)X+\mathcal{B}^\top_2\eta+R_2\bar{v}_2=\mathcal{B}^\top_2\eta+R_2\bar{v}_2,\quad a.e.,\, \mathbb{P}\mbox{-}a.s.
\end{aligned}
\end{equation}
Then
\begin{equation}
\bar{v}_2=-R^{-1}_2\mathcal{B}^\top_2\eta,\quad a.e.,\, \mathbb{P}\mbox{-}a.s.
\end{equation}
Inserting the above into (\ref{eta}), we achieve
\begin{equation}\nonumber
d\eta=-\big[(\mathcal{A}^\top+P\mathcal{F}_1)\eta+\mathcal{C}^\top\zeta-(\mathcal{F}_2+P\mathcal{B}_2)R^{-1}_2\mathcal{B}^\top_2\eta \big]ds+\zeta dW
\end{equation}
and $\eta(T)=0$. It is easy to see that $(0,0)$ is the adapted solution to the above BSDE, thus (\ref{close-loop leader}) holds. In this case,
\begin{equation}\label{y=px,z=pcx}
Y(\cdot)=P(\cdot)X(\cdot),\quad Z(\cdot)=P(\cdot)\mathcal{C}(\cdot)X(\cdot),\qquad a.e.,\, \mathbb{P}\mbox{-}a.s.
\end{equation}
and
\begin{equation}\label{u2=Theta2X}
\bar{u}_2(\cdot)=\bar{\Theta}_2(\cdot)\bar{x}(\cdot)+\bar{\tilde{\Theta}}_2(\cdot)p^2(\cdot)+\bar{v}_2(\cdot)\equiv\boldsymbol{\Theta}_2(\cdot)X(\cdot),\qquad a.e.,\, \mathbb{P}\mbox{-}a.s.
\end{equation}
Taking (\ref{close-loop leader}), (\ref{Theta2}) and (\ref{y=px,z=pcx}) into the equation of $X(\cdot)$ in (\ref{High dimension optimal system simple}), we obtain (\ref{nonanticipating X}).

From (\ref{leader riccati}), we can see that
\begin{equation}\label{cross-coupled Riccati}
\begin{cases}
\dot{P}_1+\hat{A}^\top P_1+P_1\hat{A}-P_1B_2R^{-1}_2B^\top_2P_1+C^\top P_1C+P_2\hat{F}^\top_1P_1+P_1\hat{F}_1P^\top_2+Q_2=0,\\
\dot{P}_2+\hat{A}^\top P_2+P_2\hat{A}+C^\top P_2C+P_2\hat{F}^\top_1P_2+P_1\hat{F}_1P_4-P_1B_2R^{-1}_2(B^\top_2P_2+\hat{F}^\top_2)=0,\\
\dot{P}_4+\hat{A}^\top P_4+P_4\hat{A}+C^\top P_4C+P_4\hat{F}^\top_1P_2+P^\top_2\hat{F}_1P_4-(B^\top_2P_2+\hat{F}^\top_2)^\top R^{-1}_2(B^\top_2P_2+\hat{F}^\top_2)=0,\\
P_1(T)=G_2,\quad P_2(T)=0,\quad P_4(T)=0.
\end{cases}
\end{equation}
We can see that this is a cross-coupled Riccati equation system, which solvability is difficult and we will not consider here in this paper. Finally, by It\^o's formula we have
\begin{equation}\nonumber
\begin{aligned}
&\hat{J}_2(x;\bar{\Theta}_2(\cdot)\bar{x}(\cdot)+\bar{\tilde{\Theta}}_2(\cdot)p^2(\cdot)+\bar{v}_2(\cdot))\\
&=\langle Y(0),X(0)\rangle+\mathbb{E}\int_{0}^{T} \big\langle (R_2\boldsymbol{\Theta}_2+\mathcal{F}^\top_2)X+\mathcal{B}^\top_2Y+R_2\bar{v}_2, \boldsymbol{\Theta}_2X+\bar{v}_2 \big\rangle ds\\
                                                                      &=\big\langle P_1(0)x,x\big\rangle.
\end{aligned}
\end{equation}
The proof is complete. $\qquad\Box$

\begin{Remark}
We point that here that, due to some technical reason, we could not prove the sufficiency of the above theorem, as \cite{SY2014,SY2019}. In our opinion, the completion-of-square method for the problem of the leader is invalid, because its state process triple satisfies an FBSDE, other than an SDE. How to overcome this difficulty to achieve the characterization the closed-loop solvability of Problem (SLQ)$_l$ is still open.
\end{Remark}

Finally, noting that closed-loop optimal control $\bar{u}_2(\cdot)$ of the leader has a nonanticipating representation (\ref{u2=Theta2X}) with the ``state" $X(\cdot)\equiv\left(\begin{matrix} \bar{x}(\cdot)\\ p^2(\cdot)\end{matrix}\right)$ being the solution to (\ref{nonanticipating X}). In the meanwhile, for the follower, the closed-loop optimal control $\bar{u}_1(\cdot)$ can also be represented in a nonanticipating way. In fact, from (\ref{follower closed-loop optimal}), we have
\begin{equation}
\begin{aligned}
\bar{u}_1[\bar{u}_2](\cdot)&=\bar{\Theta}_1[\bar{u}_2](\cdot)\bar{x}(\cdot)+\bar{v}_1[\bar{u}_2](\cdot)\\
                           &=-R^{-1}_1 \left(\begin{matrix} B^\top_1P^1 & 0\end{matrix}\right)X(\cdot)- R^{-1}_1\left(\begin{matrix} 0 & B^\top_1\end{matrix}\right)Y(\cdot)\\
                           &=-R^{-1}_1\Big[ \left(\begin{matrix} B^\top_1P^1 & 0\end{matrix}\right)+ \left(\begin{matrix} 0 & B^\top_1\end{matrix}\right)P \Big]X(\cdot).
\end{aligned}
\end{equation}

\section{Concluding remarks}

In this paper, we have investigated the closed-loop solution for a special LQ Stackelberg stochastic differential game. The notion of the closed-loop solvability is introduced, which require to be independent of the initial state. The follower's problem is solved first, and his closed-loop optimal strategy is characterized by a Riccati equation, together with an adapted solution to a linear BSDE. Then the necessary conditions of the existence of the leader's nonanticipating closed-loop optimal strategy is obtained via a system of cross-coupled Riccati equations. How to obtain the sufficiency in Theorem \ref{thm4.2} remains open. The solvability and numerical method of the cross-coupled Riccati equation system (\ref{cross-coupled Riccati}) is interesting and challenging. We will extend these results to control-dependent diffusions, random coefficients and mean-field case in the future.

\end{document}